\documentclass[12pt]{amsart}

\usepackage{mathrsfs}
\usepackage{amsmath,amssymb,enumerate,tikz}
\usepackage{amsfonts,amssymb}
\usepackage{booktabs}
\usepackage[colorlinks,linkcolor=red,anchorcolor=blue,
citecolor=green,CJKbookmarks=True]{hyperref}
\usepackage[T1]{fontenc}
\usepackage{csquotes}

\usepackage{imakeidx}
\makeindex

\usepackage{soul, xcolor, todonotes, easyReview}

\theoremstyle{plain}
\newtheorem{theorem}{Theorem}[section]
\newtheorem{lemma}[theorem]{Lemma}

\newtheorem{proposition}[theorem]{Proposition}

\theoremstyle{definition}
\newtheorem{definition}[theorem]{Definition}
\newtheorem{remark}[theorem]{Remark}
\newtheorem{question}{Question}

\newcommand{\tec}{Teichm\"uller }

\newcommand{\dth}{d_{Th}}
\renewcommand{\stretch}{\mathbf{stretch}}

\newcommand{\ML}{\mathcal{ML}}
\newcommand{\MF}{\mathcal{MF}}
\newcommand{\PML}{\mathcal{PML}}

\newcommand{\PMF}{\mathcal{PMF}}
\newcommand{\Hopf}{\mathrm{Hopf}}

\newcommand{\C}{\mathcal{C}}
\newcommand{\T}{\mathcal{T}}
\newcommand{\R}{\mathbb R}

\renewcommand{\H}{\mathbb H}

\title{The geometry of the Thurston metric: a survey}
\author{Huiping Pan \and Weixu Su}

\address{Huiping Pan\\
School of mathematics, South China  University of Technology, 510641, Guangzhou, China}
\email{panhp@scut.edu.cn}

\address{Weixu Su\\
	School of Mathematics, Sun Yat-sen University, 510275, Guangzhou, China}
\email{suwx9@mail.sysu.edu.cn}

\begin{document}
\maketitle
\begin{abstract}
  This chapter is a survey about the Thurston metric on the Teichm\"uller space. 
  The central issue is the constructions of  extremal Lipschitz maps between hyperbolic surfaces.  We review several constructions, including  the original work of Thurston.
 Coarse geometry and isometry rigidity of the Thurston metric, relation between the Thurston metric and the Thurston compactification  are discussed. 
Some recent  generalizations and developments of the Thurston metric are sketched.

\end{abstract}

Mathematical classification (2010)

\tableofcontents


\begingroup
\let\clearpage\relax

\section{Introduction}

This survey aims to give a brief exposition of recent results on the Thurston metric. This asymmetric metric, defined by Thurston in his  1986 preprint \cite{Thurston1998}, has a lot of beautiful geometric and analytic properties. Thurston's original investigation, just like his many other masterpieces, are based on first principles. 
Until now,  that preprint \cite{Thurston1998} together with 
Chapter 8 and Chapter 9 of \cite{thurston2022geometry}, are excellent references for learning hyperbolic geometry.

Aiming for a geometric understanding of the Teichm\"uller space from the point of view of hyperbolic surfaces and Lipschitz maps,  in parallel with the point of view of Riemann surfaces and quasiconformal maps, Thurston \cite{Thurston1998} considered the following problem:
\begin{quotation}
	Given any two hyperbolic surfaces $X$ and $Y$, what is the least possible value  of the global Lipschitz constant
	\begin{equation*}
		L(\phi):=\sup_{p\neq q }\frac{d(\phi(p)),\phi(q)}{d(p,q)}
	\end{equation*}
	for a homeomorphism $\phi:X\to Y$ in a given homotopy class?
\end{quotation}
Thurston gave a geometric construction of extremal Lipschitz homeomorphisms realizing the extremal Lipschitz constant $\mathrm{Lip}(X,Y)$ in the homotopy class of the identity map, and proved that the extremal Lipschitz constant is the same as  another geometric quantity:
\begin{equation*}
	K(X,Y):=\sup_{\alpha}\frac{\ell_\alpha(Y)}{\ell_\alpha(X)}
\end{equation*}
where $\ell_\alpha(\cdot)$ is the length of the geodesic representative of $\alpha$  on the corresponding surface and  $\alpha$ ranges over all nontrivial homotopy classes of simple closed curves.

For any two points $X$ and $Y$ in the \tec space, the \emph{Thurston metric} is defined as:
\begin{equation*}
	\dth(X,Y):=\log \mathrm{Lip}(X,Y).
\end{equation*}
Among others, Thurston proved that this metric is Finslerian, asymmetric, geodesic but not uniquely geodesic \cite{Thurston1998}.

Through the efforts of many mathematicians, we now have a fairly clear picture of this metric:  the horofunction boundary, the isometry group, the coarse geometry,  the analytic construction of extremal Lipschitz maps  and geodesic rays, to name a few. Due to its geometric inspiration,  the Thurston metric has been generalized to deformation spaces of other geometric structures and provides new perspectives there.

In this survey about the Thurston metric,
we restrict our attention to recent research.
We will look more carefully at several constructions of 
the extremal Lipschitz maps. Relevant developments will be briefly sketched or summarized without proofs.
Other surveys on the Thurston metrics include \cite{PapadopoulosTheret2007,PapadopoulosSu2015,Su2016}. 
In particular, the article \cite{PapadopoulosTheret2007}
contains foundational background material, both on the Teichm\"uller metric and the Thurston metric. 
As for prerequisites, the reader is expected to be familiar with 
hyperbolic geometry and Teichm\"uller theory. We recommend 
\cite{FLP2012, hubbard2016teichmuller}. 

 Here is the outline of this paper.  In \S \ref{sec:background}, we review some of the standard facts on \tec space, including the Thurston norm. In \S \ref{sec:construction}, we provide various constructions of extremal Lipschitz maps and geodesic rays of the Thurston metric. In \S \ref{sec:harmonic:stretch}, we deal with  a special type of geodesics, called \textit{harmonic stretch lines}, and discuss two versions of the Thurston (co-)geodesic flows. In \S \ref{sec:boundary}, we describe  the Thurston boundary and its identification with the horofunction boundary. In \S \ref{sec:rigidity}, we illustrate the isometry rigidity of the Thurston metric, both globally and infinitesimally. 
In  \S  \ref{sec:coarse} we discuss the coarse geometry of the Thurston metric.
In \S \ref{sec:counting}, we are concerned with the lattice counting problem. In \S \ref{sec:shearing}, we discuss the shearing coordinates. In \S \ref{sec:generalizations}, we sketch several generalizations of the Thurston metric.

\textbf{Acknowledgements.}  We would like to thank Athanase Papadopoulos for the opportunity of writing this article, for his patience during the preparation of the manuscript, and for suggesting several improvements on the exposition.
	


\section{Backgroud}
\label{sec:background}
\subsection{\tec space}
Let $S$ be an orientable closed surface of genus at least two. The \tec space\index{Teichm\"uller space} $T(S)$ is the space of equivalence classes of complex structures on $S$, where two complex structures $X$ and $Y$ are said to be \emph{equivalent} if there exists a conformal map $X\to Y$ which is homotopic to the identity map.  By the Uniformization Theorem,   the \tec space $T(S)$ is also the space of equivalence classes of hyperbolic structures on $S$, where two hyperbolic structures $X$ and $Y$ are said to be \emph{equivalent} if there exists an isometry $X\to Y$ which is homotopic to the identity map.  For simplicity, we  denote the equivalence class of the complex/hyperbolic structure   $X$  by $X$ itself.

\subsection{Geodesic laminations and measured geodesic laminations}
Let $X\in\T(S)$ be a hyperbolic surface. A \emph{geodesic lamination}\index{geodesic lamination} $\lambda$ on $X$ is a closed subset of $X$ which can be decomposed as  a disjoint union of  simple geodesics. Typical examples are simple closed geodesics. A \emph{multicurve}\index{multicurve} is a disjoint union of  distinct simple closed geodesics.

	Let $\lambda$ be a geodesic lamination. Then
	\begin{itemize}
		\item $\lambda$ is called \emph{maximal}\index{geodesic lamination!maximal geodesic lamination} if it cuts $X$ into ideal triangles, i.e. $X\backslash\lambda$ is a union of ideal triangles;
		\item $\lambda$ is called \emph{chain-recurrent}\index{geodesic lamination!chain-recurrent geodesic lamination} if it is the Hausdorff limit of multicurves.
	\end{itemize}
A chain-recurrent geodesic lamination is said to be \emph{maximal}\index{geodesic lamination!maximal chain-recurrent geodesic lamination} if it is not contained in any other chain-recurrent geodesic lamination.
 
 One can equip the set of geodesic laminations on $X$ with the Hausdorff metric and the corresponding Hausdorff topology. With this topology, the set of chain-recurrent geodesic laminations is a \emph{closed} subset of the space of geodesic laminations \cite[Proposition 6.2]{Thurston1998}. 
 
 Notice that the space of geodesic laminations is independent of the choice of hyperbolic metrics in the following sense. 
 Let $X'\in\T(S)$ be another hyperbolic surface.  
There is a natural one-to-one correspondence between the space of (resp. chain-recurrent) geodesic laminations on $X$ and the space of (resp. chain-recurrent) geodesic laminations on $X'$.  

Let $\lambda$ be a geodesic lamination.  A \emph{transverse measure}\index{transverse measure} $\mu$ on $\lambda$ is an assignment of a Radon measure $d\mu_\tau$ on  each arc  $\tau$ transverse to $\lambda$, subject to the following conditions:
\begin{itemize}
	\item if the arc $\tau'$ is a subarc of $\tau$ then $\left.d\mu_{\tau'}=d\mu_\tau\right|_{\tau'}$,
	\item if two arcs $\tau$ and $\tau$ are homotopic through a family of arcs transverse to $\lambda$, then the homotopy sends $d\mu_{\tau'}$ to $d\mu_\tau$.
\end{itemize}
A geodesic lamination $\lambda$ with a transverse measure $d\mu$  is called a \emph{measured geodesic lamination}\index{geodesic lamination!measured geodesic lamination}.  Typical examples of measured laminations are \emph{weigthed simple closed curves} and \emph{weighted multicurves}. Usually, we will simply denote a measured geodesic lamination $(\lambda,d\mu)$ by $\lambda$. For any simple closed curve $\gamma$, the \emph{intersection number}\index{intersection number} $i(\lambda,\gamma)$ is defined as:
\begin{equation*}
	i(\lambda,\gamma):=\inf \int_{\gamma'} d\mu
\end{equation*}
where $\gamma'$ ranges over all simple closed curves homopotic to $\gamma$. 

As we mentioned above, the space of geodesic laminations is independent of the choice of the hyperbolic metrics. So is the space of measured geodesic laminations. Let $\ML(S)$ be the space of measured geodesic laminations, equipped with the weak-$*$ topology: a sequence of measured laminations $\lambda_n$ converges to $\lambda$ if $i(\lambda_n,\gamma)\to i(\lambda,\gamma)$ for every simple closed curve $\alpha$. With this topology, $\ML(S)$ is homeomorphic to $\R^{6g-6+2n}$, where $g$ and $n$ are respectively the genus and the number of punctures of $S$ (see for instance \cite{PennerHarer1992}). Let  $\PML(S):=\ML(S)/_{\R^+}$ be the space of projective classes of measured laminations.  The set of weighted simple closed curves is dense in both $\ML(S)$ and $\PML(S)$.

\subsection{The Thurston metric and the maximally stretched lamination}\label{subsec:thurston:metric}
For any two hyperbolic surfaces $X$ and $Y$ in $\T(S)$, the extremal Lipschitz constant from $X$ to $Y$ is defined as
\begin{equation*}
	\mathrm{Lip}(X,Y):=\inf_f \{L(f)\}
\end{equation*}
where $L(f)$ is the Lipschitz constant of $f$ and where $f$ ranges over all Lipschitz homeomorhpsims from $X$ to $Y$ in the homotopy class of the identity map.  The Thurston metric\index{Thurston metric}  $d_{Th}$ is defined as
\begin{equation}\label{eq:def:th:lip}
	d_{Th}(X,Y):=\log \mathrm{Lip}(X,Y)  
\end{equation}
 Thurston \cite[Theorem 8.5]{Thurston1998} gave another characterization of this metric in terms of simple closed curves:
\begin{equation}\label{eq:def:th:length}
	d_{Th}(X,Y)=\log \sup_{\alpha} \frac{\ell_\alpha(Y)}{\ell_\alpha(X)}
\end{equation}
where $\ell_\alpha(\cdot)$ represents the length of the (unique) geodesic representative of $\alpha$ and where $\alpha$ ranges over all homotopically nontrivial simple closed curves. Since weighted simple closed curves are dense in $\ML(S)$,   one can rewrite the above identity as:
\begin{equation*}
	d_{Th}(X,Y)=\log \max_{\alpha\in\ML(S)} \frac{\ell_\alpha(Y)}{\ell_\alpha(X)}.
\end{equation*}
In general, the measured laminations realizing the maximal ratio of length are not unique. For instance, if  $\lambda$ is a non-uniquely ergodic measured lamination which realizes the maximal ratio of length from $X$ to $Y$, then any measured lamination with the same support as $\lambda$ also realizes the maximal ratio from $X$ to $Y$. Nevertheless, Thurston introduced the notion of  \emph{maximal ratio-maximizing chain-recurrent geodesic lamination} (the \emph{maximally stretched lamination} for short) from $X$ to $Y$, which is unique. 

Given $X$ and $Y$ in $\T(S)$, the \emph{maximally stretched lamination}\index{maximally stretched lamination} $\Lambda(X,Y)$ from $X$ to $Y$ is the largest chain recurrent geodesic lamination $\lambda$  with with following property: there exists an $\exp(d_{Th}(X,Y))$-Lipschitz map, homotopic to the identity map,  from a neighbourhood of $\lambda$ on $X$ to a neighbourhood of $\lambda$ on $Y$ which  takes leaves of $\lambda$ on $X$ to corresponding leaves of $\lambda$ on $Y$ by multiplying arc length by  a factor of $\exp(d_{Th}(X,Y))$. It is also  the union of all chain recurrent geodesic laminations with the aforementioned property \cite[Theorem 8.2]{Thurston1998}.  Equivalently,  the maximally stretched lamination is also  the union of   chain-recurrent laminations to which the restriction of every  Lipschitz map from $X$ to $Y$ with the extremal Lipschitz constant, which is homotopic to the identity,  takes leaves of $\lambda$ on $X$ to corresponding leaves of $\lambda$ on $Y$ by multiplying arclength by  a factor of $\exp(d_{Th}(X,Y))$.  (For the equivalence between these two descriptions, we refer to \cite[Section 9]{GueritaudKassel2017}, see also \cite[Section 5]{DaskalopoulosUhlenbeck2022} for a related discussion.)

Regarding the behavior of maximally stretched laminations, we have:

\begin{theorem}[\cite{Thurston1998}, Theorem 8.4]
Let $X,Y\in\T(S)$ be any two distinct hyperbolic surfaces. If $X_i$ and $Y_i$ are sequences of hyperbolic surfaces converging to $X$ and $Y$ respectively, then $\Lambda(X,Y)$ contains any lamination in the limit set of  $\Lambda(X_i,Y_i)$ in the Hausdorff topology.
\end{theorem}

Any Lipschitz map $f$ from $X$ to $Y$  realizing the extremal Lipschitz constant $\mathrm{Lip}(X,Y)$ is called an \emph{extremal Lipschitz map}\index{extremal Lipschitz map}. If moreover, $f$ maximally stretches exactly along the maximally stretched lamination $\Lambda(X,Y)$, then it is called an \emph{optimal Lipschitz map}\index{optimal Lipschitz map} from $X$ to $Y$. 

Before closing this subsection, let us mention the following result:
\begin{theorem}[\cite{Thurston1998}, Theorem 8.5] Let $X,Y\in  \T(S)$.
\begin{enumerate}
	\item   The extremal Lipschitz constant $\mathrm{Lip}(X,Y)$ can always be realized by a homeomorphism. 
	\item There exists a Thurston geodesic from $X$ to $Y$.
\end{enumerate}
\end{theorem}

\begin{remark}
	(i) For various constructions of \enquote{extremal Lipschitz maps} and Thurston geodesic rays, see Section \ref{sec:construction}. (ii) In general, there may be more than one geodesic from $X$ to $Y$. The union of all such geodesics is called the \emph{envelope} from $X$ to $Y$. For more information about envelopes, we refer to \cite{DLRT2020}, \cite{Bar-NatanThesis}, and  \cite{PanWolf2023}. 
\end{remark}

\subsection{The Thurston norm}
Let $X\in\T(S)$. The \emph{Thurston norm}\index{Thurston norm} on  the tangent space $T_X\T(S)$ is defined as:
 \begin{equation}\label{eq:thurston:norm}
     \|v\|_{\mathrm{Th}}:=\sup _{\lambda \in \PML(S)}
     (d_X \log\ell_\lambda)[v], \qquad \forall~ v\in T_X \T(S_{g,n}),
   \end{equation}
   where $\PML(S)=\ML(S)/_{\R_{>0}}$ is the space of projective measured laminations and $\ML(S)$ is the space of measured laminations on $S$.  The Thurston norm is the infinitesimal norm of the Thurston metric (\cite[Page 20]{Thurston1998}, see also \cite[Theorem 2.3]{PapadopoulosSu2015}).


    The unit sphere in the cotangent bundle has a very nice description as follows.
 \begin{theorem}[\cite{Thurston1998}, Theorem 5.1]\label{thm:convex:model}
   For any hyperbolic surface $X$ of finite type, the map
   \begin{equation*}
     \begin{array}{cccc}
       d\log\ell: & \PML(S) & \longrightarrow & T^*_X\T(S)\\
        & [\mu] & \longmapsto & d_X \log \ell_{\mu}
     \end{array}
   \end{equation*}
   embeds $\PML(S)$ as the boundary of a convex neighbourhood of the origin.  This convex neighbourhood is dual to the unit ball $\{v\in T_X\T(S):\|v\|_{\mathrm{Th}}\leq1\}$.
 \end{theorem}
 
  There are flat places on the unit sphere $T^1_X\T(S)$ of $T_X\T(S)$.  A \textit{facet} is a maximal flat portion of $T^1_X\T(S)$  which has maximum possible dimension so that it has interior.

  \begin{theorem}[\cite{Thurston1998},Theorem 10.1]
    \label{thm:mostly:facets}
   There is a bijection between the set of facets on the unit sphere  $T^1_X\T(S)$ and the set of simple closed curves. In other words, every facet is contained in a plane $d_X \log\ell_\alpha=1$ for some simple closed curve $\alpha$.
  \end{theorem}
  For each simple closed curve $\alpha$, let
  $$ F_X(\alpha):=\{v\in T^1_X\T(S):~(d_X \log\ell_\alpha)(v)=1 \} $$
  be the corresponding facet  obtained in Theorem \ref{thm:mostly:facets}. 

\begin{lemma}[\cite{Pan2020}, Lemma 4.2; \cite{HOP2021}, Corollary 6.13]\label{lem:disjointness}
    Let $\alpha,\beta$ be two simple closed curves. Then
    \begin{equation*}
     i(\alpha,\beta)=0 \iff \partial \mathrm{F}_X(\alpha) \cap \partial \mathrm{F}_X(\beta)\neq \emptyset.
    \end{equation*}
  \end{lemma}
  
  The    unit sphere of the Thurston norm in the tangent space encodes not only the intersection patterns of simple closed curves, but also their hyperbolic lengths.
 In \cite{DLRT2020},  Dumas--Lenzhen--Rafi--Tao    proved that  one can recognize lengths and intersection numbers of curves on $X$ from the lengths of facets in the  unit sphere in $T_X\T(S_{1,1})$, based on the observation that  extremal points of any facet in the unit sphere in $T_X\T(S_{1,1})$  are exactly those vectors tangent to the Thurston stretch lines  directed by the two canonical completions of the  corresponding simple closed curve.  This idea is further exploited and extended to surfaces of higher complexity by Huang--Ohshika--Papadopoulos \cite{HOP2021}.   Huang--Ohshika--Papadopoulos  first observed that  stretch vectors with respect to maximal chain-recurrent laminations necessarily arise as the Hausdorff limit of a sequence of shrinking faces of the unit spheres in $T_X\T(S)$, and then proceeded to prove that one can recognize lengths of simple closed curves  on $X$ using such stretch vectors.
  
    \begin{remark}
  (i) In his thesis,  Bar-Natan \cite{Bar-NatanThesis} proved that the set of tangent vectors to the stretch lines corresponding to completions of maximal
chain-recurrent laminations is precisely the set of extreme points in the unit sphere of $\T_X\T(S)$. (ii)  For more information about the unit sphere in $T_X\T(S)$, we refer to \cite{HOP2021}. (iii)
  	For the comparison of the Thurston norm with the \tec norm, the Weil--Petersson norm, and the earthquake norm, we refer to \cite{HOPP2023}.  \end{remark}

\section{Extremal Lipschitz maps and  geodesic rays}
\label{sec:construction}

In this section, we shall sketch various constructions of extremal Lipschitz maps between hyperbolic surfaces, and the corresponding geodesic rays in the Teichm\"uller space.

\subsection{Deforming polygons via explicit homeomorphisms}

One way of constructing extremal Lipschitz maps between surfaces is to consider a pair of orthogonal measured (partial) foliations and to control the  effect of the map on the leaves of these foliations.  

Following Thurston's idea, one can deform hyperbolic surfaces by  cutting the surface into polygons, deforming each polygon individually, and then gluing the deformed polygons to  get  new surfaces. One of the advantages of this consctruction is that the involved deformations/homeomorphisms can be explicitly written down.  In the following, we shall introduce this idea with more details, starting with Thurston's original construction. 

\subsubsection{Thurston's construction.}   Let $X$ be a complete hyperbolic surface of finite area. Let $\lambda$ be a \emph{maximal geodesic lamination} on $X$, that is, the complementary region $X\backslash\lambda$ is the union of \emph{ideal triangles}. Suppose further that  $\lambda$ admits at least one leaf which does not go to a cusp on both sides.    Let $\Delta$ be an arbitrary ideal triangle in $X\backslash\lambda$. Notice that each corner of $\Delta$ can be foliated by horocycles. Extend these foliations until they fill all but the region in the center bounded by three horocycles.  Let $F$ be the extended partial foliation, called the \emph{horocycle foliation} of $\Delta$. Let $G$ be the partial foliation whose leaves are geodesics rays orthogonal to the leaves of $F$. Consider the  $e^t$-Lipschitz (self)homeomorphism $f_{t}:\Delta\to\Delta$ which fixes the central region, and  maps a horocycle which has distance $r$ from the central region to the horocycle with distance $e^t r$ linearly with respect to arc length on the horocycles, see Figure \ref{fig:thurston}.  In particular, $f_t$ sends leaves of $F$ to (other)  leaves of $F$ and keeps each leaf of $G$ invariant (as a set). Thurston proved that the collection of pullback metrics on $X\backslash\lambda$  by all such $\{f_{t}\}$  extends to a new hyperbolic metric on $X$, denoted by $\stretch(X,\lambda,t)$.

\begin{figure}
	\begin{tikzpicture}[scale=3]
		\draw[red] (0,1) arc (180:240:1.73cm)
		arc (60:120:1.73cm) arc (300:360:1.73cm); 
      \draw[fill=gray,line width=1pt]  (0.23,0.14) arc (120:180:0.46cm) arc (0:60:0.46cm) arc (240:300:0.46cm);
      \draw[line width=1pt]
      (0.26,0.08).. controls (0.22,0.07) and (0.08,-0.18)..(0. 08, -0.27)  
       (0.32,0).. controls (0.22,-0.1) and (0.18,-0.18)..(0. 18, -0.28)  
      (0.37,-0.07).. controls (0.32,-0.1) and (0.25,-0.2)..(0. 27, -0.29)  
      (0.42,-0.13).. controls (0.35,-0.2) and (0.35,-0.25)..(0. 35, -0.3) ;
      
        \draw [line width=1pt] 
         (-0.26,0.08).. controls (-0.22,0.07) and (-0.08,-0.18)..(-0. 08, -0.27) 
         (-0.32,0).. controls (-0.22,-0.1) and (-0.18,-0.18)..(-0. 18, -0.28) 
         (-0.37,-0.07).. controls (-0.32,-0.1) and (-0.25,-0.2)..(-0. 27, -0.29)  
         (-0.42,-0.13).. controls (-0.35,-0.2) and (-0.35,-0.25)..(-0. 35, -0.3) ;
        
        \draw[line width=1pt] 
         (-0.16,0.28).. controls (-0.02,0.24) and (0.06,0.25)..(0. 15, 0.29) 
        (-0.2,0.2).. controls (-0.15,0.15) and (0.15,0.15)..(0. 2, 0.2) 
        (-0.1,0.43).. controls (0,0.4) and (0.04,0.43)..(0.09, 0.45) 
         (-0.12,0.35).. controls (-0.1,0.33) and (0.1,0.33)..(0.12, 0.35) ;
        
        \draw[dashed, red] (0,1)..controls (0.01,0.5) and (0.05, 0.3)..(0.1, 0.09)  (0,1)..controls (-0.01,0.5) and (-0.05, 0.3)..(-0.1, 0.09);
        
        \draw[dashed, red] (-0.87,-0.5)..controls (-0.44,-0.24) and (-0.28, -0.15)..(-0.12, 0.05)  (-0.87,-0.5)..controls (-0.43,-0.26) and (-0.25, -0.16)..(-0.03, -0.13);

          \draw[dashed,red] (0.87,-0.5)..controls (0.43,-0.26) and (0.23, -0.19)..(0.03, -0.13)  (0.87,-0.5)..controls (0.44,-0.24) and (0.28, -0.1)..(0.13, 0.04);
          \draw (0.18,-0.4)..controls(0.19,-0.35) and (0.34,-0.4)..(0.35,-0.33) (0.18,-0.4)..controls(0.17,-0.35) and (0.02,-0.4)..(0,-0.33);
          \draw (0.18,-0.4) node[below]{$r$};
    		\draw[red] (0+2,1) arc (180:240:1.73cm)
		arc (60:120:1.73cm) arc (300:360:1.73cm); 
      \draw[fill=gray, line width=1pt]  (0.23+2,0.14) arc (120:180:0.46cm) arc (0:60:0.46cm) arc (240:300:0.46cm);
      \draw[line width=1pt] (0.62+2,-0.32).. controls (0.62+2,-0.33) and (0.60+2,-0.32)..(0. 59+2, -0.37)  
      (0.52+2,-0.24).. controls (0.5+2,-0.26) and (0.48+2,-0.28)..(0. 48+2, -0.33)  
      (0.42+2,-0.13).. controls (0.35+2,-0.2) and (0.35+2,-0.25)..(0. 35+2, -0.3) 
      (0.31+2,-0.01).. controls (0.25+2,-0.03) and (0.18+2,-0.2)..(0. 18+2, -0.27) ;
      
        \draw[line width=1pt] (-0.62+2,-0.32).. controls (-0.62+2,-0.33) and (-0.60+2,-0.32)..(-0. 59+2, -0.37)
           (-0.52+2,-0.24).. controls (-0.5+2,-0.26) and (-0.48+2,-0.28)..(-0. 48+2, -0.33)  
              (-0.42+2,-0.13).. controls (-0.35+2,-0.2) and (-0.35+2,-0.25)..(-0. 35+2, -0.3)
                (-0.31+2,-0.01).. controls (-0.25+2,-0.03) and (-0.18+2,-0.2)..(-0. 18+2, -0.27) ;
        
        \draw [line width=1pt](-0.03+2,0.70).. controls (-0.01+2,0.69) and (0.01+2,0.69)..(0.03+2, 0.7)  
        (-0.05+2,0.56).. controls (-0.04+2,0.54) and (0.04+2,0.54)..(0.05+2, 0.56) 
         (-0.1+2,0.43).. controls (0+2,0.4) and (0.04+2,0.43)..(0.09+2, 0.45) 
           (-0.15+2,0.28).. controls (-0.1+2,0.23) and (0.1+2,0.23)..(0.15+2, 0.28);
        
        \draw[dashed, red] (0+2,1)..controls (0.01+2,0.5) and (0.05+2, 0.3)..(0.1+2, 0.09)  (0+2,1)..controls (-0.01+2,0.5) and (-0.05+2, 0.3)..(-0.1+2, 0.09);
        
        \draw[dashed, red] (-0.87+2,-0.5)..controls (-0.44+2,-0.24) and (-0.28+2, -0.15)..(-0.12+2, 0.05)  (-0.87+2,-0.5)..controls (-0.43+2,-0.26) and (-0.25+2, -0.16)..(-0.03+2, -0.13);

          \draw[dashed,red] (0.87+2,-0.5)..controls (0.43+2,-0.26) and (0.23+2, -0.19)..(0.03+2, -0.13)  (0.87+2,-0.5)..controls (0.44+2,-0.24) and (0.28+2, -0.1)..(0.13+2, 0.04);
           \draw (0.18+2+0.11,-0.45)..controls(0.19+2+0.11,-0.38) and (0.34+2+0.23,-0.47)..(0.35+2+0.23,-0.4) (0.18+2+0.11,-0.45)..controls(0.17+2+0.18,-0.38) and (0.02+2,-0.4)..(0+2,-0.32);
          \draw (0.18+2+0.12,-0.45) node[below]{$e^tr$};
 \draw[->, thick] (0.8,0.2)--(1.3,0.2);
 \draw (1.1,0.2)node[above]{$f_{t}$};
	\end{tikzpicture}
	\caption{In both pictures, the black arcs are leaves of the horocycle foliation $F$ while  the red dashed arcs are leaves of $G$.}
	\label{fig:thurston}
\end{figure}
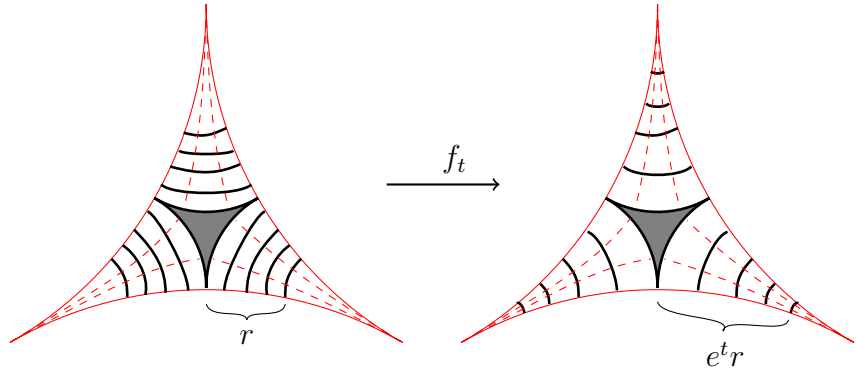

\begin{theorem}[\cite{Thurston1998}]\label{thm:thurston:stretch}
	For any complete hyperbolic surface $X$ of finite area, for any maximal geodesic lamination $\lambda$ not all of whose leaves go to a cusp on both sides, there exists a new hyperbolic surface $\stretch(X,\lambda,t)$ depending analytically on $t>0$, such that
	\begin{enumerate}
		\item the identity map $X\to \stretch(X,\lambda,t)$ is Lipschitz with Lipschitz constant $e^t$, and
		\item the identity map exactly expands arc length of $\lambda$ by the constant factor $e^t$.
	\end{enumerate}
\end{theorem}

The family $\{\stretch(X,\lambda,t):t>0\}$ is called the  \emph{Thurston stretch ray}\index{Thurston stretch ray} directed by $\lambda$. It is a geodesic under the Thurston metric.  For the proof of Theorem \ref{thm:thurston:stretch}, we refer to \cite[Section 4]{Thurston1998} and \cite[Chapter 2, Section 3.5]{PapadopoulosTheret2007}. 

\begin{remark}
	Thurston stretch lines are rare in the sense that at any point $X\in\T(S)$, the set of unit vectors tangent  to Thurston stretch lines through $X$ is of   Hausdorff dimension zero \cite[Theorem 10.5]{Thurston1998}. However, for any pair of distinct points in $\T(S)$, there exists a geodesic from the first point to the second point which a concatenation of finitely many stretch segments \cite[Theorem 8.5]{Thurston1998}
\end{remark}

\begin{remark}
Calderon and Farre \cite[Lemma 15.11 and Proposition 15.12]{CalderonFarre2021} generalized Thurston's construction to the case where $X\backslash\lambda$ is a union of regular ideal polygons.
\end{remark}

\subsubsection{The construction of Papadapoulos-Th\'eret}  A right-angled hyperbolic hexagon is said to be \emph{symmetric}\index{symmetric right-angled hyperbolic hexagon} if it has three nonadjacent edges of equal length.    Using \emph{symmetric right-angled hyperbolic hexagons}, Papadopoulos and Th\'eret  \cite{PapadopoulosTheret2012} constructed the \emph{first examples} of Thurston geodesics which are also geodesics in the reversed direction (up to reparametrization).  Their construction  starts by constructing extremal Lipschitz homeomorphisms between \emph{symmetric right-angled hyperbolic hexagons} with controlled Lipschitz constant. By doubling the hexagons, they then obtain extremal Lipschitz homeomorphisms between \emph{symmetric hyperbolic pairs of pants}, that is, hyperbolic pairs of pants which have three geodesic boundary components of equal length. By gluing pairs of pants along their boundary components without twist and by combining the resulting Lipschitz maps between pairs of pants, they obtain new stretch lines in the Teichm\"uller space of any hyperbolic surface of finite type which are also geodesics (up to reparametrization) in the opposite direction.

The Lipschitz homeomorphisms between symmetric right-angled hyperbolic hexagons is constructed as follows.  
Consider a symmetric right-angled hyperbolic hexagon $H$. There are two triples of pairwise non-consecutive edges, one of which is called \emph{long} with common length $L$ while the other is called \emph{short} with common length $l$. Similarly as in the case of ideal triangles, one can equip $H$ with two orthogonal partial measured foliations $(F,G)$ defined as follows. The   leaves of $F$ are the loci of equidistant points from the short edges. Those leaves foliate $H$ except a triangle-shaped central region.  The transverse measure on $F$ is induced by the arclength of the long edges of $H$.  The leaves of $G$ are geodesics orthogonal to the leaves of $F$ and the transverse measure of $G$ is induced by the arclength of short edges of $H$ (see Figure \ref{fig:PT}.) 

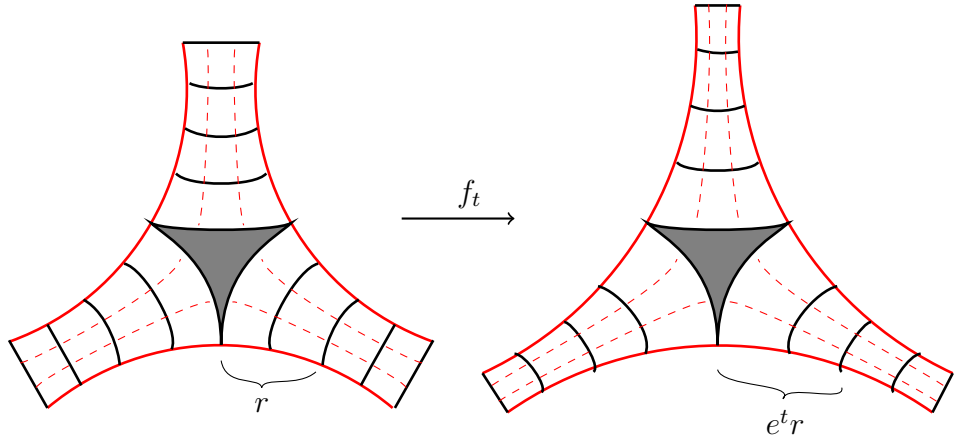
\begin{figure}
	\begin{tikzpicture}[scale=3]
		\draw[line width=1pt, red](0.17,0.98) arc (170:250: 1.19) (-0.17,0.98) arc (10:-70: 1.19) (-0.77,-0.64) arc(130:50:1.19); 
		\draw[line width=1pt] (0.17,0.98)--(-0.17,0.98) (-0.94,-0.34)--(-0.77,-0.64) (0.77,-0.64)--(0.94,-0.34);
		\draw[fill=gray, line width=1pt](0.31,0.18)..controls (0.25,0.12) and (0,0)..(0,-0.36)--
		(0,-0.36) ..controls  (-0,0) and (-0.25,0.12)..(-0.31,0.18)--(-0.31,0.18)..controls (-0.24,0.14) and (0.24,0.14)..(0.31,0.18)--cycle;
		\draw [line width=1pt](-0.2,0.4)..controls (-0.15,0.34) and (0.15,0.34)..(0.2,0.4)
		 (-0.16,0.6)..controls (-0.08,0.55) and (0.08,0.55)..(0.16,0.6)
		  (-0.14,0.8)..controls (-0.08,0.77) and (0.08,0.77)..(0.14,0.8);
		\draw [line width=1pt]
		(0.43,0)..controls (0.34,-0.02) and (0.18,-0.3)..(0.22,-0.38) 
		(0.61,-0.16)..controls (0.5,-0.22) and (0.45,-0.4)..(0.45,-0.45) (0.77,-0.27)..controls (0.6,-0.57) and (0.63,-0.52)..(0.62,-0.54);
		
				\draw [line width=1pt]
		(-0.43,0)..controls (-0.34,-0.02) and (-0.18,-0.3)..(-0.22,-0.38) 
		(-0.61,-0.16)..controls (-0.5,-0.22) and (-0.45,-0.4)..(-0.45,-0.45) (-0.77,-0.27)..controls (-0.6,-0.57) and (-0.63,-0.52)..(-0.62,-0.54);
		
		\draw[dashed, red] (-0.89,-0.44)..controls (-0.85,-0.44) and (-0.23,-0.12)..(-0.18,0.02) (-0.82,-0.55)..controls (-0.8,-0.5) and (-0.13,-0.16)..(-0.05,-0.17);

		\draw[dashed, red] (0.89,-0.44)..controls (0.85,-0.44) and (0.23,-0.12)..(0.18,0.02) (0.82,-0.55)..controls (0.8,-0.5) and (0.13,-0.16)..(0.05,-0.17);
		
		\draw[dashed, red] 
		(0.06,0.95)..controls (0.05,0.84) and (0.07,0.25)..(0.1,0.17)
		(-0.06,0.95)..controls (-0.05,0.84) and (-0.07,0.25)..(-0.1,0.17);
		
		   \draw (0.18,-0.4-0.15)..controls(0.19,-0.35-0.15) and (0.34,-0.4-0.15)..(0.42,-0.33-0.15) (0.18,-0.4-0.15)..controls(0.17,-0.35-0.15) and (0.02,-0.4-0.15)..(0,-0.33-0.1);
          \draw (0.18,-0.55) node[below]{$r$};

			\draw[line width=1pt, red](0.1+2.2,1.145) arc (175:245: 1.642) (-0.1+2.2,1.145) arc (5:-65: 1.642) (0.952+2.2,-0.66) arc (55:125: 1.642); 
		 
		\draw[line width=1pt] (0.1+2.2,1.145)-- (-0.1+2.2,1.145)
		(0.952+2.2,-0.66)--(1.042+2.2,-0.486)
		(-0.93+2.2,-0.66)--(-1.042+2.2,-0.486);
		\draw[fill=gray, line width=1pt](0.314+2.2,0.181)..controls (0.25+2.2,0.12) and (0+2.2,0)..(0+2.2,-0.362)--
		(0+2.2,-0.362) ..controls  (-0+2.2,0) and (-0.25+2.2,0.12)..(-0.31+2.2,0.18)--(-0.31+2.2,0.18)..controls (-0.24+2.2,0.14) and (0.24+2.2,0.14)..(0.31+2.2,0.18)--cycle;

  \draw [line width=1pt](-0.18+2.2,0.45)..controls (-0.15+2.2,0.4) and (0.15+2.2,0.4)..(0.18+2.2,0.45)
		 (-0.12+2.2,0.7)..controls (-0.05+2.2,0.67) and (0.05+2.2,0.67)..(0.12+2.2,0.7)
		  (-0.09+2.2,0.95)..controls (-0.02+2.2,0.93) and (0.02+2.2,0.93)..(0.09+2.2,0.94);
		
		\draw [line width=1pt]
		(0.53+2.2,-0.1)..controls (0.44+2.2,-0.12) and (0.28+2.2,-0.3)..(0.32+2.2,-0.4) 
		(0.68+2.2,-0.26)..controls (0.65+2.2,-0.25) and (0.51+2.2,-0.45)..(0.55+2.2,-0.48) (0.9+2.2,-0.4)..controls (0.85+2.2,-0.4) and (0.75+2.2,-0.55)..(0.78+2.2,-0.57);
		
				\draw [line width=1pt]
		(-0.53+2.2,-0.1)..controls (-0.44+2.2,-0.12) and (-0.28+2.2,-0.3)..(-0.32+2.2,-0.4) 
		(-0.68+2.2,-0.26)..controls (-0.65+2.2,-0.25) and (-0.51+2.2,-0.45)..(-0.55+2.2,-0.48) (-0.9+2.2,-0.4)..controls (-0.85+2.2,-0.4) and (-0.75+2.2,-0.55)..(-0.78+2.2,-0.57);
		
		\draw[dashed, red] (-1+2.2,-0.54)..controls (-0.85+2.2,-0.44) and (-0.23+2.2,-0.12)..(-0.18+2.2,0.02) (-0.97+2.2,-0.61)..controls (-0.8+2.2,-0.5) and (-0.13+2.2,-0.16)..(-0.05+2.2,-0.17);

		\draw[dashed, red] (1+2.2,-0.54)..controls (0.85+2.2,-0.44) and (0.23+2.2,-0.12)..(0.18+2.2,0.02) (0.97+2.2,-0.61)..controls (0.8+2.2,-0.5) and (0.13+2.2,-0.16)..(0.05+2.2,-0.17);
		
		\draw[dashed, red] 
		(0.04+2.2,1.15)..controls (0.02+2.2,0.84) and (0.07+2.2,0.25)..(0.1+2.2,0.17)
		(-0.04+2.2,1.15)..controls (-0.02+2.2,0.84) and (-0.07+2.2,0.25)..(-0.1+2.2,0.17);

		 \draw (0.18+2+0.2+0.11,-0.45-0.15)..controls(0.19+2+0.11+0.2,-0.38-0.15) and (0.34+2+0.23+0.2,-0.47-0.15)..(0.35+2+0.2+0.2,-0.4-0.15) (0.18+2+0.11+0.2,-0.45-0.15)..controls(0.17+2+0.18+0.2,-0.38-0.15) and (0.02+2+0.2,-0.4-0.15)..(0+2+0.2,-0.32-0.15);
          \draw (0.18+2+0.32,-0.6) node[below]{$e^tr$};
 \draw[->, thick] (0.8,0.2)--(1.3,0.2);
 \draw (1.1,0.2)node[above]{$f_{t}$};
		\end{tikzpicture} 
		\caption{Each of the hexagons has three long sides of equal length (red) and three short sides of equal length (black). The black arcs are leaves of $F$, which are loci of equidistant points from the short edges.  The red dashed arcs are leaves of $G$, which are geodesics orthogonal to the leaves of $F$}
		\label{fig:PT}
\end{figure}

Now Let $H_t$ be a family of symmetric right-angled hyperbolic hexagons with long edges being of length $t>0$.  For any $t>1$, Papadopoulos and Th\'eret  constructed a homeomorphism $f_t:H\to H_t$ such that 
\begin{itemize}
	\item $f_t$ sends the central region on $H$ to the central region on $H_t$,
	\item $f_t$ sends affinely  the leaves of $F$ at distance $r$ from the central region on $H$  to the leaves at distance $tr$ from the central region on $H_t$,
	\item $f_t$ sends leaves of $G$ on $H$ to the corresponding leaves of $G$ on $H_t$ by multiplying arclength by the factor $t$,
	\item $f_t$ is a $t$-Lipschitz homeomorphism.
\end{itemize}
Doubling the symmetric right-angled hexagons along short edges and then gluing the resulting pairs of pants without twist, Papadopoulos and Th\'eret constructed a family of hyperbolic surfaces $\{X_t\}_{t>0}$ which is a Thurston geodesic in the forward direction. Notice that the set of short edges of hexagons are glued together to form a pants decomposition of $X_t$ while the set of long edges are glued together to form several closed geodesics orthogonal to the aforementioned pants decomposition.   Interchanging  the role of short edges and long edges, we see that $\{X_t\}_{t>0}$ is also a Thurston geodesic in the backward direction. 

\subsubsection{The construction of Papadopoulos-Yamada} Starting with an \emph{arbitrary} right-angled hyperbolic hexagon, Papadopoulos and Yamada constructed a one-parameter family of right-angled hexagons with a Lipschitz map between any two elements in this family, realizing the smallest Lipschitz constant in the homotopy class of this map relative to the boundary 
\cite{PapadopoulosYamada2017}.  Let $H$ be an \emph{arbitrary} right-angled hexagon. Notice that there are two triples of non-consecutive edges on $H$, one of which called \emph{long} while the other is called \emph{short}. Let $\alpha_1,\alpha_2,\alpha_3$ be the three {long edges} of $H$.  Considering the equidistant curves to the long edges of $H$ and the corresponding orthogonal geodesic arcs,   Papadopoulos and Yamada constructed a pair of orthogonal measured (partial) foliations $(F,G)$ on $H$, which extends the construction in \cite{Papadopoulos2012}. (The construction here is more involved depending on whether the lengths of short edges satisfy the triangle inequality or not, see \cite[Section 2]{PapadopoulosYamada2017}). Let $H_t$ be a family of right-angled hyperbolic hexagons such that 
\begin{equation*}
	\frac{\cosh \ell_{H_t}(\alpha_i)}{\cosh \ell_{H}(\alpha_i)}=t.
\end{equation*}
Let  
\begin{equation*}
	k_i(t):=\frac{\ell_{H_t}(\alpha_i)}{\ell_H(\alpha_i)}, \qquad k(t)=\max\{k_i(t)\}.
\end{equation*}
 For any $t>1$, Papadopoulos and Yamada  constructed a homeomorphism $f_t:H\to H_t$ such that 
\begin{itemize}
	\item $f_t$ sends the central region on $H$ to the central region on $H_t$,
	\item $f_t$ linearly sends the leaves of $F$ on $H$  to the corresponding leaves   on $H_t$, 
	\item $f_t$ (nonlinearly) sends  leaves of $G$ on $H$ to the corresponding leaves of $G$ on $H_t$,
	\item $f_t$ is a $k(t)$-Lipschitz homeomorphism. 
\end{itemize}
Notice that each long side $\alpha_i$ is a leaf of $F$, and hence is stretched under $f$ by the factor $k_i(t)$.

 Unlike the construction of Thurston and Papadopoulos-Th\'eret, where the constructed homeomorphism \emph{linearly} stretches  the leaves of $G$ with a \emph{common factor} and \emph{linearly} shrinks   the leaves of $F$, the map here linearly stretches leaves of $F$ with \emph{varying factors} and \emph{non-linearly} shrinks  leaves of $G$. In particular, the map $f:H\to H_t$ is linear on each of the three long sides of $H$ but  \emph{not} linear on any of the three short sides of $H$. 
 
 Now we consider an arbitrary hyperbolic surface $X$ of genus $g$ with $b$ geodesic boundary components. Let $\{\gamma_i\}_{1\leq i\leq 6g-6+3n}$ be a collection of pairwise disjoint geodesic arcs orthogonal to the geodesic boundary components of $X$, i.e. $\{\gamma_i\}_{1\leq i\leq 6g-6+3n}$ is a \emph{triangulation} of $X$.  The complementary region $X\backslash(\cup_i \gamma_i)$ is a union of right-angled hexagons $\{H^i\}$.  Viewing those $\{\gamma_i\}$ as long sides of $H^i$ and applying the previously constructed Lipschitz homeomorphisms to $H^i$, we arrive at a family of hyperbolic surfaces $X_t$ and a family of Lipschitz homeomorphisms $f_t:X\to X_t$ attaining  the optimal Lipschitz constant on some of the $\{\gamma_i\}$. Doubling $X_t$ along the geodesic boundary components, we get a family of closed hyperbolic surfaces $\{X^d_t\}$  and a family of Lipschitz homeomorphisms $f_t^d:X^d\to X^d_t$ attaining  the optimal Lipschitz constant on some of the  double of $\{\gamma_i\}$.

\subsubsection{The construction of Huang-Papadopoulos} Since the optimal Lipschitz map is not as rigid outside the maximally stretched loci, one may hope to construct Thurston geodesics which keep some non-trivial subsurface invariant. This is partially done by Huang and Papadopoulos under some conditions. 

 Extending the previous constructions to \emph{ideal Saccheri quadrilaterals}, Huang and Papadopoulos \cite[Theorem 3.3 and Theorem 4.7]{HuangPapadopoulos2019} constructed,  for any  \emph{complete} hyperbolic torus $\mathbb T$ with one  boundary component and any  \emph{chain-recurrent geodesic lamination} $\lambda$ in the convex hull of $\mathbb T$ ,  a family of \emph{complete} hyperbolic tori $\{\mathbb T_t\}_{t\geq1}$ and a family of Lipschitz homeomorphisms $f_t:\mathbb T\to \mathbb T_t$ such that
\begin{itemize}
	\item the convex hull of $\mathbb T_t$ has the same boundary length as that of $\mathbb T$,
	\item $f_t$ is $t$-Lipschitz and multiplies the arclength of $\lambda$ by the factor $t$,
	\item $f_t$ is an isometry outside a compact neighbourhood of the convex hull of $\mathbb T$,
	\item if the boundary length of the convex hull of $\mathbb T$ is at most $4\mathrm{arcsinh}(1)$, $f_t$ maps the convex hull of $\mathbb T$ onto the convex hull of $\mathbb T_t$.
\end{itemize} 
In particular, for any hyperbolic torus $T$ with one geodesic boundary component of length at most $4\mathrm{arcsinh}(1)$, the construction above gives a family of hyperbolic tori $\{T_t\}_{t\geq1}$, each of which has exactly one geodesic boundary component, and a family of Lipschitz homeomorphisms $f_t:T\to T_t$ such that
\begin{itemize} 
	\item$T_t$ has the same boundary length as that of $T$, 
	\item $f_t$ is  $t$-Lipschitz,  multiplies the arclength of $\lambda$ by the factor $t$,  and maps the boundary of $T$ linearly to the boundary of $T_t$ by the factor $1$.
\end{itemize} 

Consider a  hyperbolic surface $X$ with punctures or  geodesic boundary components. If $X$ contains a hyperbolic torus $T$ with one geodesic boundary component,  whose  boundary length is at most $4\mathrm{arcsinh}(1)$, then for any chain-recurrent geodesic lamination in $T$ the above construction gives a family of hyperbolic surfaces $X_t$ and a family of Lipschitz homeomorphisms $f_t:X\to X_t$ such that 
\begin{itemize}
	\item $f_t$ is $t$-Lipschitz, and multiplies the arclength of $\lambda$;
	\item $f_t$ is an isometry on the closure of $X\backslash T$.
\end{itemize}

\subsection{Deforming crowned surfaces via harmonic diffeomorphisms}  Using harmonic maps of high energy, Pan and Wolf constructed \emph{piecewise harmonic stretch lines} and \emph{harmonic stretch lines} \cite{PanWolf2022}. Recall that a differentiable map $f:(M,\sigma |dz|^2)\to (N,\rho |dw|^2)$ between Riemannian surfaces is said to be \emph{harmonic}\index{harmonic map} if  it satisfies the \emph{Euler--Lagrange equation}:
\begin{equation*}
	f_{z\bar z}+ (\log \rho)_w f_z f_{\bar z}=0.
\end{equation*} 
If $M$ and $N$ are compact, then $f$ is harmonic if and only if it is a critical point of the energy functional 
\begin{equation*}
	E(f):=\int_M e(f) \sigma  |dz|^2
\end{equation*}  
where $e(f):=\frac{\rho(f(z))}{\sigma(z)}(|f_z|^2+|f_{\bar z}|^2)$ is the energy density of $f$.   Notice that the energy depends on the conformal structure on $M$ and the metric on $N$.  

The basic existence result of harmonic maps was established by Eells and Sampson in  \cite{EellsSampson1964} and by  Hamilton in \cite{Hamilton1975},    if the target manifold has nonpositive sectional curvature.   The uniqueness was obtained by Al'ber \cite{alber1964} and Hartman \cite{Hartman1967} if the target manifold has negative sectional curvature and if the image is not contractible to a point or a geodesic.  Moreover, Sampson \cite{Sampson1978} and Schoen-Yau \cite{SchoenYau1978} proved that any harmonic map between compact surfaces which is homotopic to a diffeomorphism is a diffeomorphism, provided that the target surface has nonpositive curvature. 

Suppose that $f$ is harmonic.
Consider the pullback  of $\rho$ {by $f$}:
$$ f^*(\rho)= \rho f_z \overline{f_{\bar{z}}}dz^2+e(f)\sigma dz d\bar{z} +\rho \overline{f_z}f_{\bar{z}} d\bar{z}^2.$$ The $(2,0)$-part of $f^*(\rho)$ is called the \emph{Hopf differential}\index{Hopf differential} of $f$.  The harmonicity of  $f$ implies that the Hopf differential of $f$ is holomorphic (see \cite{Schoen1984,Jost1984}).  If we choose the coordinate $z=x+iy$ such that the Hopf differential  $\Phi=dz^2$ and choose $\sigma$ to be the singular flat metric induced by $|\Phi|$, then $f^*(\rho)$ can be simply expressed as 
 $$ f^*\rho=(e+2)dx^2+(e-2)dy^2. $$
 By \cite[Lemma 3.2 and Lemma 3.3]{Minsky1992},  the energy density sasifies $e(z)=2+O(e^{-2r})$ where $r$ represents the distance from  $z$ to the zeros of the Hopf differential.  In particular, for points far away from the zeros of the Hopf differential, the pullback metric $f^*\rho$ is nearly  $4dx^2$. Roughly speaking, at regions far away from the zeros of the Hopf differential, the harmonic map $f$ linearly maps horizontal leaves  of the Hopf differential to hyperbolic geodesics by a factor of $2$ while contracting the  vertical leaves exponentially. This observation is the starting point of the construction of Pan and Wolf.
 
   \begin{theorem}[\cite{PanWolf2022}, Theorem 1.7]\label{thm:generalized:stretchmap}

    Let  $Y\in\T(S)$ be any closed hyperbolic surface, and let  $\lambda$ be any geodesic lamination. Then for any surjective harmonic {diffeomorphism} $f:X\to Y\setminus \lambda$ from some (possibly disconnected) punctured surface $X$,  there is a family of new hyperbolic surfaces
    \begin{equation*}
     Y_t:= \stretch(Y,\lambda,f; t)\in\T(S)
    \end{equation*}
    depending analytically on {$\{t>0\}$} such that
    \begin{enumerate}[(a)]
    \item the identity map $f_t:X\to Y_t\setminus\lambda$ is a   surjective harmonic map $f_t:X\to Y_t\setminus\lambda$ with  Hopf differential  $t\Hopf(f)$;
      \item for any $0< s<t$, the identity map ($f_t\circ f_s^{-1}$) is $\sqrt{t/s}$-Lipschitz with (pointwise) Lipschitz constant strictly less than  $\sqrt{t/s}$ in $S-\lambda$, and
       exactly expands arc length of $\lambda$ by the constant factor $\sqrt{t/s}$.
          \end{enumerate}
  \end{theorem}

A family of hyperbolic structures $\stretch(X,\lambda; \widehat{\Phi},\widehat{f}; t)$ constructed above is called a \emph{piecewise harmonic stretch  line}\index{piecewise harmonic stretch  line}.   It admits a canonical orientation coming from the orientation of the  positive real ray $\{t>0\}$. With that orientation, a piecewise harmonic stretch line is a (reparametrized) geodesic in  the Thurston metric.  Whenever we say a piecewise harmonic stretch line, we mean a directed line.

\subsection{The constructions of Gu\'eritaud-Kassel and Alessandrini-Disarlo} Using a totally different method, Gu\'eritaud and Kassel extended Thurston's construction to a geometrically finite setting and to higher dimension. 
Let $\Gamma_0$ be a discrete group. For a pair $(j,\rho)$ of representations of $\Gamma_0$ into $\mathrm{PO}(n,1)=\mathrm{Isom}(\H^n)$ with $j$ geometrically finite. 
Gu\'eritaud and Kassel investigate the set  of $(j,\rho)$-equivariant  Lipschitz maps from $\H^n$ to itself having the minimal Lipschitz constants (such maps do exists if $\rho$ is reductive \cite[Lemma 4.10]{GueritaudKassel2017}), establish the existence of  \emph{maximally stretched locus}  $E(j,\rho)$ (which is a geodesic lamination if the minimal Lipschitz constant is at least one),  and prove the existence of \emph{optimal Lipschitz maps} whose stretched locus  is exactly the maximally stretched locus $E(j,\rho)$ of the representations $(j,\rho)$ \cite[Theorem 1.3]{GueritaudKassel2017}.  This allows Gu\'eritaud and Kassel to study proper actions of discrete subgroups of $\mathrm{PO}(n,1)\times\mathrm{PO}(n,1)$ on $\mathrm{PO}(n,1)$ by left and right multipliciation \cite[Theorem 1.8 and Theorem 1.9]{GueritaudKassel2017}. 

One of the key ideas of Gu\'eritaud and Kassel is to consider the \enquote{average of Lipchitz maps}. This idea is  used by Alessandrini and Disarlo to construct extremal Lipschitz maps between hyperbolic surfaces with geodesic boundary components \cite{AlessandriniDisarlo2019}. 

We now sketch the construction of Alessandrini and Disarlo \cite{AlessandriniDisarlo2019}. Let $X$ be a finite-area hyperbolic surface with geodesic boundary components.   A geodesic lamination $\lambda$ on $X$ is called \emph{maximal} if any connected component of its complementary region $X\backslash\lambda$ is either an ideal triangle, a  right-angled quadrilateral with two consecutive ideal vertices, a right-angled pentagon with one ideal vertex, or a  right-angled hexagon. These four types of polygons are called \emph{geometric pieces}.  For any maximal geodesic lamination on $X$, Alessandrini and Disarlo constructed a family of hyperbolic surface $X^t_\lambda$ homeomorhic to $X$ and a family of generalized stretch maps $\Psi^t:X\to X^t_\lambda$ which realizes the minimal Lipschitz constant from $X$ to $X^t_\lambda$  and maximally stretches along $\lambda$. Their construction  has three steps. 

Step 1:  \emph{Generalized stretch map between geometric pieces.}\index{generalized stretch map} Combine Thurston's stretch maps and Gu\'eritaud-Kassel's idea of averaging Lipschitz maps to  construct  implicit extremal maps between geometric pieces of the same type which stretch along certain sides.  However, it is unclear whether the resulting maps between geometric pieces are  homeomorphic  or not. So one can't simply \enquote{glue} these new pieces and the corresponding stretch maps.  To overcome this issue, Alessandrini and Disarlo adapted a different idea in the remaining steps.

Step 2: \emph{Decomposition of $X$}. Let $B\subset X$ be the union of geometric pieces of $X\backslash\lambda$ which are not ideal triangles. The surface $B$ is a crowned hyperbolic surface. Let $C\subset B$ be the union of crown ends of $B$ and set $B_C= \overline{B\backslash C}\subset B$. Finally, let $X_C:=\overline{X\backslash B_C}\subset X$.

Step 3: \emph{Define $X^t_\lambda$ and the generalized stretch map $\psi^t:X\to X^t_\lambda$.} For each $t\geq 0$, Alessandrini and Disarlo  construct hyperbolic surfaces $B^t, B_C^t, X_C^t$ homeomorphic to $B, B_C, X_C$ respectively.
The new surfaces come with preferred (locally) isometric embeddings  $\xi^t: C^t \to B^t$ and $h^t:C^t\to B_C^t$. The new surface $X^t_\lambda$ is then defined to be 
\begin{equation*}
	X^t_\lambda:= (B^t \cup B_C^t)/\sim ,
\end{equation*}
where $\xi^t(z)\sim h^t(z)$ for all $z\in C^t$.  The generalized stretch map $\Psi^t:X\to X^t_\lambda$ is defined by glueing together suitable generalized stretch maps $\beta^t:B\to B^t$ and $\psi^t$ from an open dense subset of $B_C$ to $X^t_C$. (For the details of the construction, we refer to \cite[Section 6.2, 7.4, 8.1 and 8.2]{AlessandriniDisarlo2019}.

\begin{remark}
	(1) It is unclear whether the generalized stretch maps are homeomorphisms or not. (2) Based on the generalized stretch maps and the generalized stretch lines, Alessandrini and Disarlo  proved that the Teichm\"uller space of bordered hyperbolic surfaces with the \emph{arc metric} is geodesic, and provided new geodesics of  the Thurston metric in the Teichm\"uller space of closed surface. 
\end{remark}

\subsection{The construction of Daskalopoulos-Uhlenbeck}
In \cite{DaskalopoulosUhlenbeck2020,DaskalopoulosUhlenbeck2022}, Daskalopoulos and Uhlenbeck constructed extremal Lipschitz maps between manifolds using \emph{infinity harmonic maps}. Here we sketch the construction of \cite{DaskalopoulosUhlenbeck2022}, which works more generally than hyperbolic surfaces. 

Let $(M,g)$ be a compact Riemannian manifold with boundary $\partial M$ (possibly empty) and $\mathrm{dim}(M)=n$. Let $(N,h)$ be a closed Riemannian manifold of non-positive sectional curvature and let $W^{1,p}(M,N)$ denote the Soblev space of maps $f:M\to N$ such that $f$ and its weak derivative have a finite $L^p$ norm.  For any map 
 $$f:M\to N\in W^{1,p}(M,N)\cap C^0(M,N),$$
 define the $J_p$-functional of $f$ as 
 \begin{equation*}
 	J_p(f):=\int_M \mathrm{Tr}(Q(df)^p)*1,
 \end{equation*}
 where $1\leq p<\infty$ and  where $Q(df)^2=df df^T$ is a non-negative symmetric linear map mapping the  tangent space $TN$ to itself. The integrand $\mathrm{Tr}(Q(df)^p)$ is essentially  the sum of the $p$-th powers of the singular eigenvalues of $df$.  The Euler-Lagrange equations of the $J_p$ functional are
 \begin{equation*}
 	D^*Q(df)^{q-2}(df)^2=0,
 \end{equation*}
where $D=D_f$ is the pullback of the Levi-Civita connection on $f^{-1}TN$.  

If the domain manifold $(M,g)$ has nonempty boundary, there are two types minimizing problems: the Dirichlet problem and the Neumann problem.  For the Dirichlet problem we fix a continuous map $f : M \to  N$ and seek a minimizer in the homotopy class of $f$ relative to the boundary values of $f$. For the Neumann problem we only fix a homotopy class and no boundary condition at all.  Here is the existence theorem concerning the $J_p$ functional.

 \begin{theorem}[\cite{DaskalopoulosUhlenbeck2022}, Theorem 1.1] \label{thm:DU:Jp}
 	Let $(M,g)$ be a compact Riemannian manifold with boundary $\partial M$ (possibly empty) and $\mathrm{dim}(M)=n$. Let $(N,h)$ be a closed Riemannian manifold of non-positive sectional curvature. Then for each $p>N$, there exists a minimizer $f$ in $W^{1,p}(M,N)$ of the functional $J_p$ with either the Dirichlet  or the Neumann boundary conditions in a homotopy  (relative homotopy) class. Furthermore, if $\partial M\neq \emptyset$, the solution of Dirichlet problem is unique.
 \end{theorem}

 The next result is  about the construction of extremal Lipschitz maps using $J_p$ harmonic maps. Recall that a Riemannian manifold $M$ with boundary is said to be \emph{convex} if it can be isometrically embedded as a convex subset of a manifold of the same dimension without boundary. 
 
  \begin{theorem}[\cite{DaskalopoulosUhlenbeck2022}, Theorem 1.4]\label{thm:DU:extremal}
 	Let $M, N$ be Riemannian manifolds where $M$ is convex with possibly non-empty boundary and $N$ is closed and has non-positive curvature. Given a sequence $p \to \infty$, there exists a subsequence (denoted also by $p$) and a sequence of $J_p$ minimizers $f_p : M \to  N$ in the same homotopy class (either in the absolute sense or relative to the boundary depending on the context) such that
 	\begin{equation*}
 		f=\lim_{p\to\infty } f_p \in W^{1,s}(M,N)\text{ for all } s.
 	\end{equation*}
 	Furthermore, $f_p\to f$ uniformly and $f$ is an extremal Lipschitz map in the (relative) homotopy class.
 \end{theorem}

A critical point of the $J_p$ functional (with either the Dirichlet problem or the Neumann problem) is called a \emph{Schatten-von Neumann $p$-harmonic map}\index{harmonic map!Schatten-von Neumann $p$-harmonic map} or simply a $J_p$ \emph{harmonic map}\index{harmonic map!$J_p$ harmonic map} . The limit map obtained in Theorem \ref{thm:DU:extremal} is called \emph{infinity harmonic map}\index{harmonic map!infinity harmonic map}.
 About the regularity, Daskalopoulos and Uhlenbeck conjecture that for $p>2$,  any $J_p$ harmonic map between \emph{hyperbolic surfaces} is of $C^{1,\alpha}$ \cite[Conjecture 4.19]{DaskalopoulosUhlenbeck2022}.  
 
 We now return to the hyperbolic surface setting. Let $f:M\to N$ be a Lipschitz map between hyperbolic surfaces $M$ and $N$. Here $M$ and $N$ are not necessary homeomorphic.   Suppose further that the minimal Lipschitz constant in the homotopy class of $f$ is at least one.\footnote{If the minimal Lipschitz constant is less than one, then the maximally stretched loci may not be a geodesic lamination \cite[Section 9.4]{GueritaudKassel2017}. In fact, Gu\'eritaud and Kassel conjectured that in this case the maximally stretched loci is a \enquote{gramination} (contraction of \enquote{graph} and \enquote{lamination}) \cite[Conjecture 1.4]{GueritaudKassel2017}.}  In this setting, the maximally stretched locus,  i.e. the intersection of maximally stretched loci of all extremal Lipschitz maps from $M$ to $N$ in the given homotopy class, is a chain recurrent geodesic lamination.  In particular, the maximally stretched loci of  the infinity harmonic map obtained in Theorem \ref{thm:DU:extremal} contains the maximally stretched lamination.  For the inverse direction, they proved:
 \begin{theorem}[\cite{DaskalopoulosUhlenbeck2022}, Theorem 1.7]\label{thm:DU:canonical}
 	If the maximally stretched lamination $\lambda$  for a homotopy class of maps $M\to N$ consists of a finite number of simple closed geodesics, then the maximally stretched loci of the infinity harmonic  map obtained in Theorem \ref{thm:DU:extremal} is exactly $\lambda$. 
 \end{theorem}
 
 Daskalopoulos and Uhlenbeck conjectured that Theorem \ref{thm:DU:canonical} holds in general \cite[Conjecture 8.12]{DaskalopoulosUhlenbeck2022}: the maximally stretched loci of any infinity harmonic map, whose Lipschitz constant is at least one, is the maximally stretched lamination \cite[Conjecture 8.12]{DaskalopoulosUhlenbeck2022}. Furthermore, they conjectured that the limit infinity harmonic maps are not homeomorphisms \cite[Conjecture 8.5]{DaskalopoulosUhlenbeck2022}.

\subsection{Concatenations of geodesic segments} 
In the previous subsections, we discussed the constructions of various geodesic rays. Thurston constructed another type of geodesics, that is,  concatenations of Thurston stretch segments \cite[Theorem 8.5]{Thurston1998}. This construction applies to other segments as well.

\section{Harmonic stretch lines}\label{sec:harmonic:stretch}
In Section \ref{sec:construction}, we discussed various constructions of geodesic rays under the Thurston metric.  In this section, we will focus on a special type of piecewise harmonic stretch lines, namely, the \emph{harmonic stretch lines}. Throughout this section, we assume that the underlying topological surface $S$ is a closed orientable surface with genus at least two.

\subsection{Harmonic stretch lines}
 Recall that  the construction of piecewise harmonic stretch lines through $Y\in\T(S)$ depends on the choice of a geodesic lamination $\lambda$ and a surjective harmonic diffeomorphism  $f:X\to Y\backslash\lambda$ from some (possibly disconnected)  punctured surface $X$. If $X$  has several components, the harmonic map $f$ can be prescribed independently on each component. If we impose some additional \enquote{compatibility conditions} on the  harmonic maps from distinct components, we arrive at  the \emph{harmonic stretch lines}. Originally, the harmonic stretch lines are defined via harmonic map rays, of which we now recall the definition.

 Let  $X\in\T(S)$  and let $\Phi$ be a  holomorphic quadratic differential on $X$.  The \emph{harmonic map ray}\index{harmonic map!harmonic map ray} determined by $X$ and $\Phi$ is the set of hyperbolic surfaces $\{Y_t\}_{t>0}$ such that the Hopf differential of harmonic map from $X$ to $Y_t$ is $t\Phi$ \cite{Wolf1989}.

  A piecewise harmonic stretch line is said to be a \emph{harmonic stretch line}\index{harmonic stretch line} if it arises as a limit of some sequence of harmonic map rays in $\T(S)$.  The basic theorem about those harmonic stretch lines is the following.

 \begin{theorem}[\cite{PanWolf2022}, Theorem 1.8]\label{thm:unique:geodesic}
    For any two distinct hyperbolic surfaces $Y, Z\in \T(S)$, there exists a unique harmonic stretch line proceeding from $Y$ to $Z$. 
  \end{theorem}

 Actually, one can say more about the relationship between harmonic map rays and Thurston geodesic rays. To that end, we need another definition, the \emph{harmonic map dual rays}.  Let $X\in\T(S)$ and let $\lambda$ be a measured foliation on $X$. The \emph{harmonic map dual ray} determined by $X$ and $\lambda$ is the of hyperbolic surfaces such that the horizontal foliation of the Hopf differential of the harmonic map from $X_t$ to $X$ is $t\lambda$ (for the existence and uniqueness of such a $X_t$, see \cite[Theorem 3.1]{Wolf1998}). For any two hyperbolic surfaces $X',X$, we denote by $\mathbf{HR}(X',X)$ the harmonic map which starts at $X'$ and passes through $X$.
 
 \begin{theorem}[\cite{PanWolf2022}, Theorem 1.1]\label{thm:pw:hr}
 	Let $X\in \T(S)$ be a hyperbolic surface and $\lambda$ a measured lamination. Then the family of  harmonic map rays $\mathbf{HR}(X_t,X)$ converge to a Thurston geodesic locally uniformly as $X_t$ diverges along the harmonic map dual ray determined by $X$ and $\lambda$.
 \end{theorem}
 \begin{remark}
 The conclusion still holds if we replace the harmonic map dual ray by a Teichm\"uller ray \cite[Theorem 1.2]{PanWolf2022}.
 Similar results also hold for harmonic map dual rays. 
 	We fixed $X\in\T(S)$ and a holomorphic quadratic differential $\Phi$ on $X$. Suppose that  $Y_t$ diverges along the  harmonic map ray determined by $X$ and $\Phi$. Then the family of harmonic map dual rays, each of which starts at $Y_t$ and passes through $X$, converges to the  Teichm\"uller geodesic determined by $X$ and $\Phi$ locally uniformly \cite[Theorem 1.4]{PanWolf2022}.
 \end{remark}
 
 If we remove the assumption that $X_t$ diverges along a harmonic map ray, then we have the following sub-convergence result:
 \begin{theorem}[\cite{PanWolf2022}, Theorem 1.3]
 	For any fixed $X\in\T(S)$, let 
 	 $X_n\in\T(S)$ be an arbitrary divergent sequence.  with $X_n\to\infty$. Then  the  sequence of harmonic map rays $\mathbf{HR}(X_n,X)$ contains a subsequence which converges to some Thurston geodesic locally uniformly.
 \end{theorem}
 
  Next, we extend the existence/uniqueness theory of harmonic stretch lines to rays whose terminal point is a projective measured lamination, representing a point on the boundary of the Thurston compactification of the \tec space $\T(S)$.  

   \begin{theorem}[\cite{PanWolf2022}, Theorem 1.11]\label{thm:exponentialmap}
   For any $Y\in \T(S)$ and  any $[\eta]\in\PML(S)$, there exists a unique harmonic stretch ray  starting at $Y$, which converges to $[\eta]\in\PML(S)$ in the Thurston compactification.
   
   Moreover, these rays foliate $\T(S)$ if we fix $Y$ and let $[\eta]$ vary in $\PML(S)$, or if we fix $[\eta]$ and let $Y$ vary in $\T(S)$. 
 \end{theorem}

 Finally, let us mention the following continuity property, where $\mathrm{HSL}(X,X')$ represents the harmonic stretch line proceeding from $X\in\T(S)$ to $X'\in\T(S)\cup\PML(S)$
 \begin{proposition}[\cite{PanWolf2022}, Proposition 12.14 and Proposition 13.8] 
 \label{prop:hsl:continuity}Let $Y$ and $Z$ be two distinct points in $\T(S)$.  Let $[\eta]\in\PML(S)$.
 	\begin{enumerate}[(i)]
 		\item Suppose that $Y_n\to Y$ and $Z_n\to Z$ in $\T(S)$, then the sequence of harmonic stretch lines $\mathrm{HSL}(Y_n,Z_n)$ locally uniformly converges  to the harmonic stretch line $\mathrm{HSL}(Y,Z)$.
 		\item  Suppose that $Y_n\to Y$ in $\T(S)$ and $[\eta_n]\to [\eta]$ in $\PML(S)$, then the sequence of harmonic stretch lines  $\mathrm{HSL}(Y_n,[\eta]_n)$ locally uniformly converges  to the harmonic stretch line   $\mathrm{HSL}(Y,[\eta])$.
 	\end{enumerate}
 \end{proposition}

  \begin{remark}
 Using the identification between $\ML(S)$ and $T^*_Y\T(S)$, Theorem \ref{thm:exponentialmap} induces a bijection from  $T^*_Y\T(S)$ to $\T(S)$ which sends rays through the origin in $T^*_Y\T(S)$ to harmonic stretch rays in $\T(S)$. (By Proposition \ref{prop:hsl:continuity} such a bijection is a homeomorphism.) 
 	A natural question is to improve Theorem  \ref{thm:exponentialmap} to an exponential map from $T_Y\T(S)$ to $\T(S)$. Since the Thurston norm is not strictly convex, this does not follow directly from Theorem \ref{thm:exponentialmap}.  The existence part of the exponential map is clear, see \cite[Remark 1.12]{PanWolf2022}. However, the uniqueness part is unknown yet. 
 \end{remark}

 \subsection{Two versions of the geodesic flow\index{Thurston geodesic flow} for the Thurston metric}
 
 Theorem \ref{thm:exponentialmap} allows us to define  the Thurston geodesic flow  $$\psi_t: \T(S)\times \PML(S)\longrightarrow \T(S)\times\PML(S)$$
 such that the orbit through $(Y,[\eta])\in \T(S)\times \PML(S)$ is the harmonic stretch line determined by $Y$ and $[\eta]$ via Theorem \ref{thm:exponentialmap}. Moreover, every harmonic stretch line appears as a (forward) orbit.

 There is another version of the Thurston geodesic flow  
 $$\phi_t:\T(S)\times\PML(S) \to \T(S)\times\PML(S)$$
 such that the orbit through $(Y,[\eta])\in \T(S)\times\PML(S)$ is the stretch line obtained from Theorem \ref{thm:pw:hr}. 

For a further discussion about  those flows, we refer to \cite{PanWolf2022}.

\section{Thurston boundary}

\label{sec:boundary}

The Teichm\"uller space
$\T(S)$ parametrizes hyperbolic structures on the topological surface $S$. A sequence of hyperbolic surfaces $X_n\in \T(S)$ which 
tends to infinity has the property that some closed curve $\gamma$ on
$S$ becomes
infinitely long, that is, the hyperbolic length
$\ell_\gamma(X_n)\to \infty$.

The length of long curves can be measured quantitatively by looking at their transverse measures with respect to measured foliations. 
 This observation of Thurston \cite{Thurston1986} resulted in a natural compactification of $\T(S)$, which is called the Thurston compactification. 

Recall that a \emph{measured foliation}\index{measured foliation} $(F, \mu)$ is a foliation $F$ on $S$ equipped with
 an invariant transverse measure
$\mu$. We require that the singularities of $F$ are similar to the singularities of holomorphic quadratic differentials
(or meromorphic quadratic differentials with at most simple pole at the punctures, when $S$ is a punctured surface).

For simplicity, we denote $(F, \mu)$ by $\mu$.
For any simple closed curve $\gamma$, the \emph{intersection number}\index{intersection number} $i(\gamma, \mu)$ is defined  by

$$i(\gamma,\mu)=\inf_{\gamma'}\int_{\gamma'} d\mu,$$
where $\gamma'$ is taken over all simple closed curves homotopic to $\gamma$.
Two measured foliations $\mu_1, \mu_2$ are said to be equivalent if $i(\gamma,\mu_1)=i(\gamma, \mu_2)$
for all simple closed curves $\gamma$ on $S$. Let $\MF$  be the set of equivalence classes of measured
foliations on $S$. The set of projective equivalence classes of measured foliations is denoted by $\PMF(S)$.

As shown by Thurston (see \cite{FLP2012}), $\PMF(S)$ is topologically a sphere of dimension $6g-7+2n$.
Moreover, $\PMF(S)$ is the Thurston boundary of $\T(S)$. By definition, a sequence of hyperbolic surfaces
$X_n\in \T(S)$ converges to the projective class of $\mu$ if there is some $c_n>0$ such that
$$c_n \ell_\gamma(X_n) \to i(\gamma,\mu)$$
for all simple closed curves  $\gamma$ on $S$. To glue $\T(S)$ together with $\PMF(S)$, we can embed both $\T(S)$ and $\PMF(S)$ into $\mathrm{P}\mathbb{R}^\infty$.

There are many applications of Thurston's theory of measured foliations. Below, we will limit ourselves to illustrating some of the results related to the Thurston metric. 

The following result, first proved by  Papadopoulos \cite[Theorem 5.1]{Papadopoulos1991},
implies there is direct relation between the Thurston metric and the Thurston boundary:

\begin{theorem}
Let $\lambda$ be a maximal geodesic lamination on $X\in \T(S)$ with $\mu=F_\lambda(X)$ being the transverse horocycle foliation.
Then the corresponding Thurston stretch ray converges to the projective class $[\mu]\in \PMF(S)$.
\end{theorem}

\begin{remark}
The problem of the convergence of anti-stretch lines (stretch lines in backward direction) towards the Thurston boundary was studied by Th\'eret \cite{Theret2007}. If the anti-stretch line is directed by a complete geodesic lamination $\lambda$ which is the completion of a uniquely ergodic measured lamination $\mu$, then it converges to the projective class of $\mu$.	For general stretch lines, the problem of the convergence or non-convergence remains open. 
\end{remark}


It is true that every geodesic ray of the Thurston metric has a unique limit in the Thurston boundary.
This is a  corollary of the following result by Walsh \cite[Theorem 3.6]{Walsh2014}.

\begin{theorem}\label{thm:horo}
The horofunction boundary\index{horofunction boundary} of $\T(S)$, endowed with the Thurston metric, is homeomorphic to the Thurston boundary.
\end{theorem}
The horofunction compactification and the corresponding boundary of a proper metric
space are introduced by Gromov \cite{Gromov1981}.

\begin{remark}\label{rmk:arcmetric}
An asymmetric metric, called \textit{arc metric}, is defined on the Teichm\"uller space of compact hyperbolic surfaces with boundary.
 It is an analogue of the Thurston  metric defined on the Teichm\"uller space of complete finite area hyperbolic surfaces without boundary.
 Theorem \ref{thm:horo} was generalized to the arc metric by \cite{ALPS2016}.
Geodesics for the arc metric,  which are analogues of Thurston's stretch lines, were constructed by \cite{AlessandriniDisarlo2019}. By \cite[Theorem 14.1]{PanWolf2022}, the arc metric coincides with the Lipschitz metric, which is defined using Lipschitz maps similarly as \eqref{eq:def:th:lip}. 
\end{remark}

Let us briefly mention the work of Bonahon \cite{Bonahon1988} on geodesic currents. Let us endow the surface $S$ with a hyperbolic metric. For simplicity, we assume that $S$ is a closed hyperbolic surface. 

A \textit{geodesic current}\index{geodesic current} is a locally finite Radon (positive) measure on the space of geodesics on $S$. Alternatively, a geodesic current is a finite Radon measure on the unique tangent bundle of $S$ that is invariant under the geodesic flow. Closed geodesics or weighted sums of closed geodesics on $S$ are  examples of geodesic currents. 

Denote the space of geodesic currents on $S$ by $\mathcal{C}(S)$, endowed with the weak$^*$ topology. We can embed $\MF(S)$ into $\mathcal{C}(S)$, by representing each measured foliation as a geodesic lamination with an invariant transverse measure. 
On the other hand, every hyperbolic structure $X\in \T(S)$  corresponds to a \textit{Liouville current}\index{Liouville current}. By taking projective classes, the Teichm\"uller space embeds into $\mathrm{P}\mathcal{C}(S)$, and Bonahon \cite{Bonahon1988} proved that the  boundary of the image is exactly $\PMF(S)$.

  \section{Isometry Rigidity}
  \label{sec:rigidity}
  In this section, we will discuss the rigidity of the Thurston metric, especially the isometry rigidity.  
  
  \subsection{Global rigidity}

    We start with the global rigidity due to Walsh \cite{Walsh2014}.  Using an identification  of  the horofunction boundary of the Thurston metric with the Thurston boundary (Theorem \ref{thm:horo}), Walsh proved that
  \begin{theorem}[\cite{Walsh2014}, Theorem 7.8]\label{thm:Walsh1}
  	If $S$ is not a sphere with four or less punctures, nor a torus with two or fewer punctures, then every isometry of $(\T(S),d_{Th})$ is an element of the extended mapping class group.
  \end{theorem}
  Theorem \ref{thm:Walsh1} is an analogue of Royden's theorem for the Teichm\"uller metric.
Note that the proof of Walsh is geometric. An important step consists of reducing the problem into showing that every isometry of the Thurston metric induces an automorphism of the curve complex.  
  
  Furthermore, distinct surfaces give rise to distinct Teichm\"uller spaces:
  \begin{theorem}[\cite{Walsh2014}, Theorem 7.8]
  	Let $S_{g,n}$ and $S_{g',n'}$ be  orientable surfaces of finite type with negative Euler characteristics.  Assume that $\{(g,n),(g',n')\}$ is not one of the three sets:
  	\begin{equation*}
  		\{(1,1),(0,4)\},~\{(1,2),(0,5)\},~\{(2,0),(0,6)\}.
  	\end{equation*} 
  	If $(g,n)\neq (g',n')$, then $(\T(S),d_{Th})$ and $(\T(S_{g',n'}),d_{Th})$ are not isometric.  
  \end{theorem}

  \subsection{Infinitesimal rigidity} The infinitesimal rigidity of the Thurston metric was first studied by Dumas--Lenzhen--Rafi--Tao \cite{DLRT2020} for once-punctured tori, then independently by Pan \cite{Pan2020} and Huang-Ohshika-Papadopoulos \cite{HOP2021} for surfaces of higher complexity.
\begin{theorem}[\cite{DLRT2020, Pan2020, HOP2021}]
\label{thm:inf:rigidity}
Let $S$ be a surface of finite type with negative characteristic.
  Let $X,Y\in\T(S)$. Then there exists an isometry of normed vector spaces 
 \begin{equation*}
(T_X\T(S),\|\bullet\|_{\mathrm{Th}})
  \to(T_Y\T(S),\|\bullet\|_{\mathrm{Th}})
 \end{equation*}
  if and only  $X$ and $Y$ are in the same extended mapping class group orbit.
\end{theorem}

We can also use the Thurston  norm to distinguish the topology of the underlying hyperbolic surfaces.
\begin{theorem}[\cite{Pan2020}]\label{thm:topology:rigidity}
	Let $S$ and $S'$ be  two surfaces of finite type with negative Euler characteristics.  Then there exists an isometry of  normed vector spaces 
 \begin{equation*}
(T_X\T(S),\|\bullet\|_{\mathrm{Th}})
  \to(T_Y\T(S'),\|\bullet\|_{\mathrm{Th}})
 \end{equation*} 
 for some $X\in\T(S)$ and $Y\in\T(S')$
  if and only  $S$ and $S'$ are homeomorphic.
  \end{theorem}
\begin{remark}
	The original statement of Theorem \ref{thm:topology:rigidity} is slightly differently from the statement here. In \cite[Theorem 1.6]{Pan2020}, the isometry is assumed to be $\R$-linear. As pointed out in \cite[Remark 1.16]{HOP2021}, this assumption is unnecessary, as isometries of the Thurston norm are also isometries of its (additive) symmetrisation, and the Mazur--Ulam theorem  ensures that any isometry of the Thurston norm must be affine, and hence there exists an isometric translate of the affine isometry which is a linear map.
\end{remark}

  Using the  infinitesimal rigidity, we have the local rigidity:
  \begin{theorem}[\cite{DLRT2020, Pan2020, HOP2021}]\label{thm:local:rigidity}
  Let $S$ be a surface of finite type with negative Euler characteristic.
  	Let $U\subset \T(S)$ be a connected open subset, considered as a metric space with the restriction of the Thurston metric.  Then any isometric embedding $(U,\dth)\to (\T(S),\dth)$ is the restriction to $U$ of an element of the extended mapping class group.
  \end{theorem}
  Intuitively, this says that the quotient of $\T(S)$ by the mapping class group is \enquote{totally unsymmetric}: each ball fits into the space isometrically in only one place. As a direct consequence of Theorem \ref{thm:local:rigidity}, this reproves Theorem \ref{thm:Walsh1} and extends it to the exceptional surfaces:
   \begin{theorem}[\cite{Walsh2014, DLRT2020,Pan2020, HOP2021}]\label{thm:global:rigidity}
   Let $S$ be a surface of finite type with negative characteristic.
  	Then every isometry of $\T(S),\dth)$ is induced by an element of the extended mapping class group. 
  \end{theorem}

  To pass from the infinitesimal rigidity to the local rigidity, one needs the following regularity result about the Thurston norm: 
  \begin{theorem}[\cite{DLRT2020}, Theorem 6.1]
  	Let $S$ be a surface of finite type with negative characteristic. Then the Thurston norm function $T\T(S)\to \R$ is locally Lipschitz.
  \end{theorem}

      \section{Coarse Geometry}\label{sec:coarse}
     In this section, we will discuss the coarse geometry of the Thurston metric.  
     \subsection{Short markings}
     A \emph{pants decomposition}\index{pants decomposition} on $S$ is a collection of mutually disjoint curves which cut $S$ into pairs of pants. A \emph{marking}\index{marking} $\mu$ on $S$ is a pants curve system $P$ with additionally a set of transverse curves $Q$ satisfying the following properties. We require each curve $\alpha\in P$ to have a unique transverse curve $\beta\in Q$  that intersects $\alpha$  minimally (once or twice) and is disjoint from all other curves in $P$. We will often say $\alpha$ and $\beta$ are dual to each other, and write $\bar \alpha=\beta$ or $\bar\beta=\alpha$.       
     
     Given $X\in\T(S)$, a \emph{short marking}\index{marking!short marking} $\mu_X$ on $X$ is a marking where the pants curve system is constructed using the algorithm that picks the shortest curve on $X$, then the second shortest disjoint from the first, and so on. Once the pants curve system is complete, the transverse curves are then chosen to be as short as possible. Note that a short marking on $X$ may not be unique, but all short markings on $X$ form a bounded set in the curve complex. Thus, we will refer to $\mu_X$ as the associated short marking on $X$.
     
     \subsection{Curve graphs}
     Given two curves $\alpha$ and $\beta$ on $S$, we define their \emph{intersection number}\index{intersection number}$i(\alpha,\beta)$ to be the minimal number of intersections between any representatives of homotopy classes of $\alpha$ and $\beta$. Notice that the minimum  of intersection numbers for any two distinct (homotopy classes of) simple closed curves  is 1 for the once-punctured torus, 2 for the four-holed sphere, and 0 for all other surfaces.

     The \emph{curve graph}\index{curve graph} $C(S)$ is defined as follows: the vertices are homotopy classes of nontrivial simple closed curves and two vertices are connected by an edge if the corresponding curves can be realized to have minimal possible intersection numbers.     We equip $C(S)$ with a metric by assigning length one to every edge.
     
     If $S$ is an annulus,   the above definition does not give  interesting objects since there is only one homotopy class of simple closed curves on any annulus. We need a different definition. To emphasize the difference, we denote by $A$  the (topological)  annulus.  By an arc on $A$ we always mean a homotopy class of a simple arc $\omega$ connecting the two boundary components of $A$ where the homotopy is taken relative to the endpoints of $\omega$. The intersection $i(\omega,\omega')$ of two arcs is the minimal number of intersections between any representatives of homotopy classes of $\omega$ and $\omega'$. The vertices of  the curve graph $C(A)$ are arcs on $A$ and the edges are pairs of arcs with zero intersection. We also equip $C(A)$ with a metric as above.

     \subsection{Subsurface projection} 
     For any subsurface $Y$ of $S$, the \emph{subsurface projection}\index{subsurface projection} 
     \begin{equation*}
     	\pi_Y:C(S)\to \mathcal{P}(C(Y))
     \end{equation*}
     from $C(S)$ to the set of subsets of the vertices of $C(Y)$ is defined as follows. 
     Suppose first that $Y$ is not an annulus. If $\alpha$ is disjoint from $Y$, then $\pi_Y(\alpha)=\emptyset$. If $\alpha$ is contained in $Y$, then $\pi_Y(\alpha)=\alpha$. If $\alpha$ intersects $Y$, then $\alpha\cap Y$ is a union of arcs.  Let $\omega$ be such an arc. Let $\delta$ and $\delta'$ be the boundary components of $Y$ intersecting $\omega$.  Let $\alpha_\omega$ be the collection of boundary components of the $\epsilon$-regular neighbourhood of $\delta\cup\omega\cup \delta'$ for sufficiently small $\epsilon$. The projection $\pi_Y(\alpha)$ is then defined as the collection of  $\alpha_\omega$ for all arcs of $\alpha\cap Y$.  From the construction we see that the diameter of $\pi_Y(\alpha)$ is at most two.
     
     We now consider the case where $S$ is an annulus $A$ with core curve $\gamma$. The Gromov compactification of the  annulus cover  of $S$ corresponding to $\gamma\in\pi_1(S)$ is well-defined and is independent of the choice of the hyperbolic metrics on $S$. For any simple closed curve $\alpha\in C(S)$, the projection $\pi_A(\alpha)$ is defined to be the lifts of $\alpha$ to the annulus cover that connects two boundary components of the annulus. Notice that any lift has  exactly two well-defined endpoints in the Gromov compactification. The diameter of $\pi_A(\alpha)$ is also at most two.
     
     The projection of multicurves is defined as the union of the projection of each simple closed curve.
     
     \subsection{Combinatorial model}
     Recall that for any $X$ and $Y$ in $\T(S)$,
     \begin{equation*}
     	\dth(X,Y)=\log \sup_{\alpha} \frac{\ell_\alpha(Y)}{\ell_\alpha(X)},
     \end{equation*}
   where $\alpha$ ranges over all essential simple closed curves on $S$.  This implies that there is a simple closed curve $\alpha$  such that $\log\frac{\ell_\alpha(Y)}{\ell_\alpha(X)}$ is a good estimate of $\dth$. Lenzhen, Rafi, and Tao proved that such a curve can be chosen within any short marking $\mu_X$ of $X$:
     \begin{theorem}[\cite{LRT2012}, Theorem E]\label{thm:combinatorial}
     There exists a constant $C$ depending on the topology of $S$ such that for any $X$ and $Y$ in $\T(S)$,
     \begin{equation*}
     	\left|\dth(X,Y)-\log \max_{\alpha\in\mu_X} \frac{\ell_\alpha(Y)}{\ell_\alpha(X)}\right|\leq C.
     \end{equation*}
     \end{theorem}
     Theorem \ref{thm:combinatorial} is a consequence of the following estimate of lengths of simple closed curves in terms of short markings.
     \begin{proposition}[\cite{LRT2012}, Proposition 3.1]\label{prop:length:combinatorial}
     	There exists a positive constant $C$ depending on the topology of $S$, such that for any $X\in\T(S)$, any simple closed curve $\gamma$, and any short marking $\mu_X$ of $X$,
     	\begin{equation*}
     C^{-1}{\sum\limits_{\alpha\in \mu_X}i(\gamma,\alpha)\ell_{\bar\alpha}(X)}\leq 	{\ell_\alpha(X)}\leq C{\sum\limits_{\alpha\in \mu_X}i(\gamma,\alpha)\ell_{\bar\alpha}(X)}	     .	\end{equation*}
     \end{proposition}
     \begin{remark}
     	(1) A similar estimate also holds for the extremal length:
     	   	\begin{equation*}
     C^{-1}{\sum\limits_{\alpha\in \mu_X}i(\gamma,\alpha)^2\mathrm{Ext}_{\bar\alpha}(X)}\leq 	{\mathrm{Ext}_\alpha(X)}\leq C{\sum\limits_{\alpha\in \mu_X}i(\gamma,\alpha)^2\mathrm{Ext}_{\bar\alpha}(X)}	.	\end{equation*}
     \end{remark}
     (2) For similar estimates using twists with respect to some pants decomposition, see \cite{Minsky1996}. 
     
 \subsection{Short curves} One way to understand the behavior of Thurston geodesics in the moduli space of Riemann surfaces is to characterize short curves along the corresponding family of hyperbolic surfaces. 
Such a characterization is known for Teichm\"uller geodesics (see \cite[Theorem 1.1]{Rafi2005} or \cite[Theorem 2.3]{LRT2012}). However, the situation is more involved for Thurston geodesics.  

 \begin{definition}
 	\begin{enumerate}
 		\item Let $K$ be a positive constant. Two points $X$ and $Y$ in $\T(S)$ are said to have \emph{$K$-bounded combinatorics}\index{bounded combinatorics} if 
 		\begin{equation*}
 			d_{C(\Sigma)}(\pi_\Sigma(\mu_X),\pi_\Sigma(\mu_Y))\leq K
 		\end{equation*}
 		for any subsurface $\Sigma\subset S$, where $\mu_X$ and $\mu_Y$ are respectively the short markings of $X$ and $Y$.
 		\item A path in $\T(S)$ is said to be $\epsilon$ \emph{cobounded}\index{cobounded} if it is contained in the $\epsilon$ thick part of $\T(S)$.
 	\end{enumerate}
 \end{definition}

\begin{theorem}[\cite{LRT2012}, Theorem A]\label{thm:cobounded}
	Let $X$ and $Y$ be two points in the $\epsilon$-thick part of $\T(S)$. If $X$ and $Y$ have $K$-bounded combinatorics, then every Thurston geodesic from $X$ to $Y$ fellow travels the \tec geodesic from $X$ to $Y$, and hence is $\epsilon'$-cobounded, where $\epsilon'$ depends on $\epsilon$ and $K$. 
\end{theorem}

This reflects a negative-curvature phenomenon of the Thurston metric. More precisely,
\begin{theorem}[\cite{LRT2012}, Theorem C and Corollary D]
Let $\mathcal G$ be a Thurston geodesic in $\T(S)$ whose endpoints have bounded combinatorics. Then the closed point projection to $\mathcal G$ is strongly contracting. In particular, any quasi-geodesic with the same endpoints as $\mathcal G$ fellow travels $\mathcal G$.
\end{theorem}

\begin{remark} 
From Theorem \ref{thm:cobounded}, we see bounded combinatorics implies $\epsilon$-cobounded. One may wonder whether the inverse still holds or not. For the \tec metric, the converse also holds. In fact, two points  $X$ and $Y$ in $\T(S)$ have bounded  combinatorics if and only if the \tec geodesic connecting them is cobounded (\cite[Theorem 1.1]{Rafi2005}, see also \cite[Theorem 2.2]{LRT2012}). However, it is not the case for the Thurston metric. There exists $\epsilon_0>0$ such that for any $\epsilon>0$, there are points $X$ and $Y$ in $\T(S)$ and a Thurston geodesic from $X$ to $Y$ that stays  in the $\epsilon_0$ thick part of $\T(S)$, whereas the associated \tec geodesic from $X$ to $Y$ does not stay in the $\epsilon$-thick part of $\T(S)$ \cite[Theorem 1.4]{LRT2015}. 
 \end{remark}
 
 Recall that the Thuston metric is not uniquely geodesic. 
 	If we remove the bounded combinatorics assumption, then the fellow traveling property  fails. In fact, for any $D>0$, there are points $X,Y,Z$ in $\T(S)$ and two Thurston geodesics $\mathcal G_1$, $\mathcal G_2$ from $X$ to $Y$ with the following properties \cite[Theorem 1.1]{LRT2015}.
 	\begin{itemize}
 		\item Geodesics $\mathcal G_1$ and $\mathcal G_2$ do not fellow travel each other; the point $Z$ lies in path $\mathcal G_1$ but is at least $D$ away from any point in $\mathcal G_2$;
 		\item The geodesic  $\mathcal G_1$ parameterized in any way in the reverse direction is not a geodesic. In fact, the point $Z$ is at least $D$ away from any point in any geodesic from  $Y$ to $X$.
 	\end{itemize} 
  On the other hand, if we consider the projection to the curve complex, which  sends a point in $X\in \T(S)$ to the set of shortest simple closed curves,  then all  Thurston geodesics with common endpoints nearly look the same:
 
 \begin{theorem}[\cite{LRT2015}, Theorem 1.2]
 	The shadow of a Thurston geodesic to the curve graph is a reparameterized quasi-geodesic. 
 \end{theorem}
 
 \begin{remark}
 	Since the curve graph is Gromov hyperbolic \cite{MasurMinsky1999}, quasi-geodesics with common endpoints fellow travel each other. Hence,  for any two points $X$ and $Y$ in $\T(S)$, the shadow to the curve graph of different Thurston geodesics form $X$ to $Y$ fellow travel each other.
 \end{remark}
 
 In \cite{DLRT2020},  the authors provided another characterization of short curves along Thurston geodesics. Before stating the result, let us introduce the notion of absolute twisting  $d_\alpha(X,Y)$ of two hyperbolic surfaces $X$ and $Y$ relative to a simple closed curve $\alpha$.  Let $\omega$ be a geodesic arc on $X$  which is orthogonal to the geodesic representative of $\alpha$ on $X$. Let $\omega'$ be a geodesic arc on $Y$  which is orthogonal to the geodesic representative of $\alpha$ on $Y$. Let $\widehat X\to X$ be the annulus cover of $X$ corresponding to $\alpha$. Let $\widehat \alpha$ be the core curve of $\widehat X$.   Let $\widehat\omega$  and $\widehat{\omega}'$ be respectively   lifts of $\omega$ and $\omega'$ intersecting $\widehat X$.  The \emph{absolute twisting}\index{absolute twisting} $d_\alpha(X,Y)$ of  $X$ and $Y$ relative  to $\alpha$ is defined to be:
 \begin{equation*}
 	d_\alpha(X,Y)= \min i(\widehat \omega, \widehat \omega')
 \end{equation*}
 where the minimum ranges over all such $\omega$ and $\omega'$ and their lifts.  (Notice that $\widehat X$ has a well-defined Gromov compatification, and $i(\widehat \omega,\widehat \omega')$ is considered the minimum intersection number in the given homotopy class relative to their endpoints.)
 
 Recall that for any pair of distinct points $X$ and $Y$ in $\T(S)$, there is a well-defined lamination, the \emph{maximally stretched lamination} $\Lambda(X,Y)$ (see Section \ref{subsec:thurston:metric}). A simple closed curve $\alpha$ is said to  \emph{interact} with the maximally stretched lamination  $\Lambda(X,Y)$ if it belongs to the lamination or intersects it essentially. 
 \begin{theorem}[\cite{DLRT2020}, Theorem 1.2]
 	There exists a  constant $\epsilon_0$ such that the following statement holds. Let $X$ and $Y$ be two points in the $\epsilon_0$-thick part of $\T(S)$ and let $\alpha$ be a simple closed curve on $S$ that interacts with $\Lambda(X,Y)$. Then the minimum length $\ell_\alpha$  of $\alpha$ along any Thurston geodesic from $X$ to $Y$ satisfies:
 	\begin{equation*}
 		C^{-1} \cdot d_\alpha(X,Y)-C'\leq \frac{1}{\ell_\alpha}\mathrm{Log} \frac{1}{\ell_\alpha}\leq C \cdot d_\alpha(X,Y)+C'
 	\end{equation*} 
 	where $C$ and $C'$ are two constants that depend only on $\epsilon_0$, and where $\mathrm{Log}(x):=\max\{1,\log (x)\}$.
 	 \end{theorem}
   For more interesting discussions about short curves, we refer to \cite[Section 3]{DLRT2020}. 
   
Based on the above characterization of short curves, Telpukhovskiy proved that  Masur's criterion does not hold in the Thurston metric:
\begin{theorem}[\cite{Telpukhovskiy2022}]
	 There are Thurston stretch paths in the  Teichm\"uller space $\T(S_{0,7})$ with minimal, filling, but not uniquely ergodic horocyclic foliation, that stay in the thick part for the whole time.
\end{theorem}
     \subsection{Length spectrum metric and  the product theorem} 
      Inspired by Minsky's product theorem for the Teichm\"uller metric, Choi and Rafi proved an analogue for the symmetrized metric $d_L$ of the Thurston metric:
      \begin{equation*}
      	d_L(X,Y):=\max\{\dth(X,Y), \dth(Y,X)\}.
      \end{equation*}        
     The metric $d_L$ is called the \emph{length spectrum metric}\index{length spectrum metric}. 

     For any positive constant $\epsilon$, consider the $\epsilon$-thick part $\mathrm{Thick}_{\epsilon}(S)$ and the $\epsilon$-thin part $\mathrm{Thin}_{\epsilon}(S)$ defined as below:
     \begin{eqnarray*}
     	\mathrm{Thick}_\epsilon(S)&:=&\{X\in\T(S): \ell_\alpha(X)\geq\epsilon \text{ for every simple closed curve } \alpha\}
     	\\
     	\mathrm{Thin}_\epsilon (S)&:=&\{X\in\T(S): \ell_\alpha(X)<\epsilon \text{ for some simple closed curve } \alpha\}. 
     \end{eqnarray*}
     For a multicurve $\Gamma$, we define the subset $\mathrm{Thin}_\epsilon(\Gamma;S)$ of the $\epsilon$-thin part as the set of  points $X\in\T(S)$ satisfying $\ell_\alpha(X)<\epsilon$ if and only if $\alpha\in\Gamma$.    Let $\T(S\backslash\Gamma)$ be the \tec space of the punctured surface homeomorphic to $S\backslash\Gamma$, let $U_i:=\{(x,y)\in \R: y>1/\epsilon\}$ be the subset of the upper half plane. Then the Fenchel--Nielsen coordinates on $\T(S)$ provide a natural homeomorphism:
     \begin{equation*}
     \Pi:	\mathrm{Thin}_\epsilon(\Gamma; S) \to \T(S\backslash\Gamma)\times U_1\times \cdots \times  U_k
     \end{equation*}

     Let $d_{L(S\backslash \Gamma)}$ be the Lipschitz metric on $\T(S\backslash\Gamma)$. 
     Let $d_{L(\gamma_i)}$ be a modified metric of the hyperbolic metric on $U_i$.  Let $$d_{L_\Gamma}:=\sup \{d_{L(S\backslash \Gamma)}, d_{L(\gamma_1)}, \cdots, d_{L(\gamma_k)}\}.$$
     Choi and Rafi proved that
     \begin{theorem}[\cite{ChoiRafi2007}, Theorem C]     	There exists a constant $c$ depending on $\epsilon$ and the topology of $S$ such that for any multicurve $\Gamma$ and for any $X$ and $Y$ in $\mathrm{Thin}_\epsilon(\Gamma;S)$,  we have
     	\begin{equation*}
     	\left|	d_L(X,Y)-d_{L_\Gamma}(\Pi(X),\Pi(Y))\right|\leq c.
     	\end{equation*}
     \end{theorem}
     
     Before closing this subsection, let us mention the comparison between the length spectrum metric and the Teichm\"uller metric. In the thick part of $\T(S)$, these two metrics are almost isometric \cite[Theorem B]{ChoiRafi2007}, whereas  in the whole Teichm\"uller space   they are not comparable \cite{Li2003,ChoiRafi2007}. On the other hand,  they are almost isometric in the moduli space \cite{LiuSu2011,LSSZ2017}.

\section{Counting lattice points}\label{sec:counting}
In this section, we will discuss the lattice counting problem about the Thurston metric. Before stating the counting result precisely, we need several definitions. 

Following Bonahon \cite{Bonahon1988},  we can embed the Teichm\"uller space $\T(S)$ as a subset of  the space $\C(S)$ of geodesic currents  in such a way that for every $Y\in\T(S)$ and for every curve $\gamma$ we have
\begin{equation*}
	i(\gamma,Y)=\ell_{Y}(\gamma).
\end{equation*}
One can also embed $\T(S)$ into the space of functions on $\C(S)$ as follows.
Let $X\in\T(S)$.  Define the function $D_X$:
\begin{eqnarray*}
	D_X: \C(S)&\longrightarrow& \R \\
	\eta&\longmapsto &\sup_{\lambda\in\ML(S)}\frac{i(\eta,\lambda)}{\ell_X(\lambda)}.
\end{eqnarray*}
It is clear that $D_X:\C(S)\to\R$ is  positive, homogeneous and continuous.  

Let $\mu_{Th}$ be the Thurston volume on $\ML(S)$. For any positive, homogeneous and continuous function $f$ on $\C(S)$, define 
\begin{equation*}
	m(f):=\mu_{Th}(\{\lambda\in\ML(S):i(\alpha,\lambda)\leq 1\}).
\end{equation*}
In particular, for any $Y\in\T(S)$,
\begin{equation*}
	m(Y):=\mu_{Th}(\{\lambda\in\ML(S):\ell_Y(\lambda)\leq 1\}),
\end{equation*}
and 
\begin{equation*}
	m(D_Y):=\mu_{Th}(\{\lambda\in\ML(S):D_Y(\lambda)\leq 1\}).
\end{equation*}
Finally, let 
\begin{equation*}
	m_S:=\int_{\mathcal{M}(S)}m(Y)d\mathrm{vol}_{wp}
\end{equation*}
where $\mathcal{M}(S)$ is the moduli space and $d\mathrm{vol}_{wp}$ represents the Weil-Petersson volume.

\begin{theorem}[\cite{RafiSouto2019}, Theorem 1.1]\label{thm:counting}
	Let $S$ be a closed orientable surface of genus at least two. Let $X$ and $Y$ be two points in $\T(S)$. Then 
	\begin{equation*}
		\lim_{R\to\infty}\frac{\#\{\phi\in\mathrm{MCG}(S):\dth(X,\phi(Y))\leq R\}}{e^{(6g-6)R}}=\frac{m(D_X)m(Y)}{m_S}.
	\end{equation*}
\end{theorem}
     \begin{remark}
     	Theorem \ref{thm:counting}  is a consequence of a more general counting result due to \cite{RafiSouto2019}, where the authors counted the mapping class group orbits of  filling currents (for more information about this type of counting problems, we refer to \cite{ErlandssonUyanik2020,ErlandssonSouto2022} and references therein).  The analogue result under the Teich\"uller metric was first proved in \cite{ABEM2012}  and improved with an error term in \cite{Herrera2021}.
     \end{remark}
      

 \section{Shearing coordinates}

\label{sec:shearing}

The \emph{shearing coordinates}\index{shearing coordinates} for $\T(S)$ was introduced by Thurston \cite{Thurston1998}.
The shearing operation is a generalization of the earthquake deformation. In Thurston's preprint, the construction is called  \emph{cataclysm coordinates}.

Let $\Delta$ be an ideal triangle in the hyperbolic plane. There is a unique circle in $\Delta$
 tangent to each edge of $\Delta$.  If $\gamma$ is an edge of $\Delta$, then the circle intersects
 $\gamma$ at some distinguished point.

 Consider two ideal triangles $\Delta_a, \Delta_b$, with disjoint interiors and a common edge $\gamma$.
The shear parameter $s_\gamma$ is defined to be the signed distance of the two distinguished points on $\gamma$.
The quantity $s_\gamma$ measures the magnitude of shift (or twist) when we glue $\Delta_a$ and $\Delta_b$. We refer to \cite{BBFS2013} for details. 

Fix a finite ideal triangulation $\lambda$ of $S$. 
For each hyperbolic surface $X\in \T(S)$, $\lambda$ is represented as a geodesic lamination
with finitely many leaves whose
 complementary components are ideal triangles.
We can associate with each edge $\gamma$ of $\lambda$ a shear parameter $s_\gamma(X)$.
Then the map
\begin{eqnarray*}
s_\lambda: \T(S) &\to& \mathbb{R}^{|\lambda|} \\
X &\mapsto& \left(s_\gamma(X)\right)_{\gamma\in\lambda}
\end{eqnarray*}
is an embedding, which is called a \emph{shearing  coordinates}\index{shearing  coordinates} of $\T(S)$.
Note that the image of $s_\lambda$ is an open convex cone of a linear subspace of the  dimension
$6g+ 2n-6$.

For example, every pants decomposition of $S$ can be completed to be a finite ideal triangulation.
Then, the corresponding Fenchel--Nielsen parameters of $\T(S)$ can be extended
(up to linear transformation) to the shearing coordinates.

\begin{remark}
	In Mirzakhani's wrok, random hyperbolic surfaces are defined by gluing random pairs of pants. It is also natural to construct random hyperbolic surfaces by gluing ideal triangles randomly. 	
\end{remark}

The above construction works for general maximal geodesic laminations.
Fix a maximal geodesic lamination $\lambda$.
Then each $X\in \T(S)$ can be associated with a horocyclic measured foliation
$F_\lambda(X)$ and a \emph{shearing cocycle}\index{shearing cocycle}. 
In brief, a shearing cocycle is similar to some transverse weights on a train track carrying $\lambda$,
and the weights measure the shifts of the shear map when one hyperbolic surface is deformed into another. 

Conversely, for any shearing cocycle, one can construct a measured foliation $F_\lambda$ in a neighborhood of $\lambda$,
whose leaves are hocycles orthogonal to $\lambda$. It turns out that, there is a unique hyperbolic structure $X$ on
$S$ such that $F_\lambda$ extends to $F_\lambda(X)$.

Let $\mathcal{MF}(\lambda)$ be the set of measured foliations that is totally transverse to $\lambda$ (see \cite[Proposition 9.4]{Thurston1998}). 
Then
\begin{theorem}[Thurston]
Let $\lambda$ be a maximal geodesic lamination on $S$. Then the shear map 
\begin{eqnarray*}
F_\lambda: \T(S) &\to& \mathcal{MF}(\lambda) \\
X &\mapsto&  F_\lambda(X)
\end{eqnarray*}
is a real analytic parametrization of $\T(S)$.
\end{theorem}

The shearing coordinates can be also interpreted as transverse cycycles,
as did by Bonahon \cite{Bonahon1996}. There are several nice properties and developments of this notion:

\begin{enumerate}
  \item The shear map pulls back the Thurston form on $\MF(S)$ to the Weil--Petersson symplectic form
  on $\T(S)$ \cite{papadopoulos1993weil, BonahonSozen}. This has been very useful when Mirzakhani \cite{mirzakhani2008ergodic} proved that the earthquake flow is measured isomorphic to the horocycle flow on 
  moduli space. 
  \item The complexified version of transverse cocycles can be used to measure the bending of pleated surfaces, and thus is used
  to extend Thurston's parametrization to the quasi-Fuchsian space \cite{Bonahon1996}. 
  Furthermore, Thurston's parametrization can be extended to Hitchin representations into $\mathrm{PSL}_n(\mathbb{R})$ \cite{BonahonDreyer2014, bonahon2017hitchin}.
   
  \item There are many related works on decorated Teichm\"uller spaces, quantum Teichm\"uller spaces and higher Teichm\"uller spaces,
  such as \cite{fock2005dual,fock2006moduli,fock2007moduli}, ect. 
  
  \item There is a construction of shearing coordinates on the universal Teichm\"uller space, see 
\cite{penner1993universal, vsaric2010circle, vsaric2022circle}
\end{enumerate}

When $\lambda$ is not a maximal geodesic lamination, the horocyclic foliation is no longer well-defined. Such a difficultly was already encountered in the case of surfaces with boundaries, 
when Alessandrini--Disarlo  \cite{AlessandriniDisarlo2019} constructed geodesics for the arc metric.

In general, Calderon--Farre \cite{CalderonFarre2021} defined orthogeodesic foliations
associated with $\lambda$. 
The orthogeodesic foliations are constructed 
on complementary subsurfaces. They showed that the restriction of the orthogeodesic foliation to each component of $S\setminus\lambda$ completely determines
the hyperbolic structure on that piece. Thus they generalized the shearing cooordinates 
of $\T(S)$.

We conclude this section with the following result of Th\'eret \cite{Theret2014}, which improves the result of \cite[Theorem 1.3]{BBFS2013}.

\begin{theorem}[\cite{Theret2014}, Theorem 1.1]
	Let $\mu$ be a measured geodesic lamination transverse to the complete geodesic lamination $\lambda$ on $S$. Then the length function
	$\ell_\mu$ is a convex function of the shear coordinates defined over $\T(S)$ associated with $\lambda$, and $\ell_\mu$  is strictly convex whenever $\mu$ intersects all leaves of $\lambda$.
\end{theorem}

In particular, lengths functions are convex along stretch lines. An interesting corollary is that the length of the horocyclic foliation is strictly decreasing on a stretch line and converges to zero as $t$ converges to infinity.

\section{Generalization of the Thurston metric}

\label{sec:generalizations}

In this section, we shall discuss some generalizations of the Thurston metric. As we mentioned earlier, there is a generalization to the setting of hyperbolic surface with geodesic boundary components, the \emph{arc metric}, see Remark \ref{rmk:arcmetric}.

Gu\'eritaud-Kassel \cite{GueritaudKassel2017} considered a pair $(j,\rho)$ of representations of a discrete group $\Gamma_0$ into $G=\operatorname{Isom}(\mathbb{H}^n)
$ such that $j$ is injective with $j(\Gamma_0)$ discrete and geometrically finite.
They showed that, when the minimal Lipschitz constant is strictly greater than $1$,
there exists a geodesic lamination that is ``maximally
stretched" by any extremal $(j,\rho)$-equivariant Lipschitz map $f:\mathbb{H}^n\to \mathbb{H}^n$. As an application,
they generalized the two-dimensional results and constructions of Thurston and extend the Thurston
metric on Teichm\"uller space  to
higher dimension.

Let us denote the minimal Lipschitz constant of $(j,\rho)$-equivariant map
by $L(j,\rho)$. A related constant is 

\begin{equation*}
	K(j,\rho)=\log \sup\limits_{\gamma} \frac{\lambda(\rho(\gamma))}{\lambda(j(\gamma))},
\end{equation*}
where $\lambda$ denotes the translation length of the element, and
the supremum is taken over all hyperbolic elements in $\gamma\in\Gamma_0$.
It is obvious that
$$K(j, \rho)\leq L(j,\rho).$$

Let $M=\mathbb{H}^n/\Gamma_0$ be a geometrically finite hyperbolic manifolds.
The Teichm\"uller space $\mathcal{T}(M)$ is defined to be the set of conjugacy classes
of geometrically finite representations of $\Gamma_0$.
For $j,\rho\in \mathcal{T}(M)$, denote the critical exponents of $j$ and $\rho$ by $\delta(j)$ and $ \delta(\rho)$, respectively. Then

\begin{theorem}[\cite{GueritaudKassel2017}]
The function
$$d_{\mathrm{Th}}(j,\rho)=\log \left(L(j,\rho)\frac{\delta(\rho)}{\delta(j)}\right)$$
defines an asymsymmetric metric on $\mathcal{T}(M)$.
\end{theorem}

We call the above metric a \textit{generalized Thurston metric}. Another version is
\begin{equation}\label{equ:K}
\overline{d}_{\mathrm{Th}}(j,\rho)=\log \left(K(j,\rho)\frac{\delta(\rho)}{\delta(j)}\right).
\end{equation} 
As pointed out by Gu\'eritaud-Kassel, the two metrics differ in general.

\begin{question}
Is the generalized Thurston metric $d_{\mathrm{Th}}$ Finsler?
How to construct geodesics of $d_{\mathrm{Th}}$?
\end{question}

Recently, Sapir \cite{Sapir2022}  extended the Thurston metric to the space of projective filling currents; Carvajales-Dai-Pozzetti-Wienhard \cite{CDPW2022}
 generalized the Thurston metric and the associated Finsler norm to Anosov representations of surface groups,
including the space of Hitchin representations. The distance formula is similar to
\eqref{equ:K}.

The Thurston metric has an analogue in Outer space, see \cite{vogtmann2015geometry} for a discussion. Just like the Thurston metric on $\T(S)$, this metric is
not symmetric.  It was used by Bestvina \cite{bestvina2010bers} to give a Bers-like proof of the train track theorem  for fully irreducible automorphisms of $F_n$. 

Regarding the extremal Lipschitz maps, Thurston \cite{Thurston1998} suggested that
\begin{quotation}
	\enquote{I currently think that a characterization of minimal stretch maps should be possible in a considerably more general context (in particular, to include some version for all Riemannian surfaces)...}
\end{quotation}
For related results along this line, we refer to \cite{PanWolf2022,DaskalopoulosUhlenbeck2020, DaskalopoulosUhlenbeck2022}.

Finally, let us mention that  it is impossible to deform hyperbolic metrics  of finite area in $\T(S)$ such that no length of any closed geodesic is increased. This is a typical example of the (marked) length spectrum rigidity problem, which is extensively studied in geometry and dynamics, see \cite{otal1990spectre,duchin2010length,guillarmou2019marked} for examples.



\bibliographystyle{plain}
\bibliography{main.bib}

\printindex
\endgroup
\end{document}